\documentclass[11pt,reqno]{amsart}

\usepackage{graphicx}

\newcommand{\K}{{\mathbf K}^\crit}

\newcommand{\sm}{\setminus}

\newcommand{\kahler}{K\"ahler }

\newcommand{\PP}{{\mathbb P}}
\newcommand{\RR}{\mathbb{R}}

\newcommand{\CC}{{\mathbb C}}

\newcommand{\NN}{{\mathbb N}}

\newcommand{\dbar}{\bar\partial}
\newcommand{\ddbar}{\partial\dbar}

\renewcommand{\H}{{\mathbf H}}
\newcommand{\half}{{\frac{1}{2}}}

\newcommand{\acal}{\mathcal{A}}

\newcommand{\fcal}{\mathcal{F}}

\newcommand{\hcal}{\mathcal{H}}

\newcommand{\lcal}{\mathcal{L}}
\newcommand{\mcal}{\mathcal{M}}

\def    \deg    {{\operatorname{deg}}}

\def    \span   {{\operatorname{span}}}

\def    \half   {{\frac{1}{2}}}

\def    \span   {{\operatorname{span}}}

 \def   \half   {{\frac{1}{2}}}

 \def    \Im     {{\operatorname{Im}}}
 \def    \Re     {{\operatorname{Re}}}

\newtheorem{maintheo}{{\sc Theorem}}
\newtheorem{maincor}{{\sc Corollary}}
\newtheorem{mainprop}{{\sc Proposition}}
\newtheorem{mainlem}{{\sc Lemma}}
\newtheorem{theo}{{\sc Theorem}}[section]

\newtheorem{lem}[theo]{{\sc Lemma}}
\newtheorem{prop}[theo]{{\sc Proposition}}
\newtheorem{definition}[theo]{{\sc Definition}}
\newtheorem{exam}[theo]{{\sc Example}}

\newenvironment{defin-no-number}{\medskip\noindent{\it Definition:\/} }{\medskip}

\newtheorem{claim}[theo]{{\sc Claim}}

\def\h#1{\hbox{#1}}

\def\o{\omega}

\def\MA{Monge-Amp\`ere }

\def\K{K\"ahler }

\def\ra{\rightarrow}

\def\th{\theta}
\def\vp{\varphi}

\def\w{\wedge}
\def\i{\sqrt{-1}}
\def\text{\textstyle}

\def\ra{\rightarrow}

\def\del{\partial}

\def\MAop{\hbox{\rm MA\hglue0.02cm}}

\def\dom{\hbox{\rm dom}}

\def\dis{\displaystyle}

\def\calM{\mcal}
\def\calA{\acal}
\def\calL{\lcal}

\def\calH{\hcal} \def\H{\hcal}
\def\calF{\fcal}

\def\uscreg{{\operatorname{reg}}}
\def\cvx{{\operatorname{cvx}}}
\def\epi{{\operatorname{epi}}}

\def\usdoublestar{u^{\star\star}_s}
\let\usd=\usdoublestar
\def\reg{{\operatorname{reg}}}
\def\sing{{\operatorname{sing}}}
\def\bfT{\hbox{\bf T}}
\def\oFS{\o_{\hbox{\small FS}}}
\def\dbz{d\bar z}
\def\bigs{\bigskip}
\def\eps{\epsilon}

\font\small=cmr5

\title[
The Cauchy problem for Monge-Amp\`ere, II
]
{
The Cauchy problem for the homogeneous
Monge-Amp\`ere equation,
II. Legendre transform
}

\author{Yanir A. Rubinstein }
\author{Steve Zelditch }

\address{Department of Mathematics, Stanford University, Stanford, CA 94305, USA}
\email{yanir@member.ams.org}

\address{Department of Mathematics, Northwestern  University, Evanston, IL 60208, USA} 
\email{zelditch@math.northwestern.edu}

\thanks{\hglue-10pt October 12, 2010. Revised October 24, 2010.}

\begin{document}

\begin{abstract}
We continue our study of the Cauchy problem for the homogeneous
(real and complex) Monge-Amp\`ere equation (HRMA/HCMA). In the prequel
\cite{RZ2} a quantum mechanical approach for solving the HCMA was
developed, and was shown to coincide with the well-known Legendre
transform approach in the case of the HRMA. In this 
article---that uses tools of 
convex analysis and can be read independently---we
prove that the candidate solution produced by these methods ceases
to solve the HRMA, even in a weak sense, as soon as it ceases
to be differentiable. 

\end{abstract}

\maketitle

\tableofcontents


\bigskip
\section{Introduction and main results}
\label{IntroductionSection}

\bigs

This article is the second in a series whose aim is to study
existence, uniqueness and regularity of solutions of the 
Cauchy problem for the homogeneous (real and complex) \MA
equation (HRMA/HCMA).
Our goal in the present article is
to show that the well-known Legendre transform method
for linearizing the HRMA
fails to solve the equation as soon as the solution of
the linear equation fails to be convex.

The Cauchy problem studied in this article corresponds to the
initial value problem (IVP) for geodesics in the space of \K metrics.
The IVP can be phrased as a Cauchy problem for
the HCMA on 
the product $S_T\times M$ of a strip,
$S_T:=[0,T]\times\RR$, and a \K manifold $M$,
and in the presence of an $(S^1)^n$ symmetry the HCMA
reduces to the HRMA on $[0,T]\times\RR^n$.

In the first part \cite{RZ2}, we constructed a certain  {\it quantum
analytic continuation potential} on any projective \K manifold $M$
and conjectured that it solved
the IVP for as long as a solution exists. 
We evaluated it
explicitly in the case of torus invariant metrics on toric
and Abelian varieties, and showed that it equals
the well-known Legendre transform potential in those cases. 
The main purpose of this sequel is to investigate the lifespan 
of this Legendre transform potential. 
An important part of the analysis is to determine rather precise
regularity estimates for this function.
Our main result (Theorem \ref{FirstMainThm})
shows that it ceases to solve the HRMA, even as a weak solution,
as soon as it develops singularities, or equivalently,
when the associated solution to the linear problem ceases to be convex. 
At the same time, we show that it does solve the equation on 
its dense regular locus.
The proof of Theorem
\ref{FirstMainThm} uses only tools of convex analysis and is
independent of the quantum techniques of \cite{RZ2}.

The Cauchy problem for  HRMA studied in
this article is given by

\begin{equation}
\label{HRMARayEq}
\left\{
\begin{array}{rrl}
\h{\rm MA}\, \psi
\!\!\!& = & \!\!\!
0, \quad\quad\;\,\mskip2mu \mbox{on} \; [0,T] \times \RR^n,
\cr\cr
\psi(0,\,\cdot\,)
\!\!\!& = &\!\!\!
\psi_0(\,\cdot\,),
\;\; \mbox{on} \;  \RR^n,
\cr\cr
\displaystyle
\frac{\partial\psi}{\partial s}(0,\,\cdot\,)
\!\!\!& = &\!\!\!
\dot\psi_0(\,\cdot\,), \;\; \mbox{on} \; \RR^n.
\end{array} \right.
\end{equation}
\smallskip

\noindent
Here, $\h{\rm MA}$ denotes the real Monge-Amp\`ere operator that 
associates a Borel measure to a convex function (see \S\S\ref{MAOpSubsection})
and equals $\det\nabla^2 f\, dx^1\w\cdots\w dx^{n+1}$ on $C^2$ functions.

It is well-known that the HRMA is linearized by the Legendre
transform.
In geometric terms, the Legendre transform is an isometry
between the space of open-orbit \K potentials and the space
of symplectic potentials. Both spaces are flat, however the latter
also has a trivial connection (see \cite{RZ1}, \S3). Therefore, a
geodesic $\{\psi_s\}$ of open-orbit \K potentials, i.e., a solution
of (\ref{HRMARayEq}) with $\psi_s=\psi(s,\,\cdot\,)$, is transformed
under the Legendre transform to a straight line of
symplectic potentials
\begin{equation}
\label{UsFirstEq}
\psi_s^\star=u_s = u_0 + s \dot{u}_0
\end{equation}
defined on the polytope $P=\overline{\Im\,\nabla\psi_s}$ corresponding to the toric \K class
(see \S\S\ref{SectionCauchySympPot} for more background, as well as
\cite{RZ1}, \S3, and \cite{RZ2}, \S\S4.2).
There exists a
certain (typically) finite time $T_\span^\cvx$, that we call {\it
the convex lifespan}, such that (\ref{HRMARayEq}) restricted to
$[0,T_\span^\cvx)\times\RR^n$ may be linearized and hence solved
via the Legendre transform. One has
\begin{equation}
\label{TspanCvxEq}
T_\span^\cvx(\psi_0,\dot\psi_0):=
\sup\{s\,:\, \psi_0^\star-s\dot\psi_0\circ(\nabla\psi_0)^{-1} \h{\ is convex on
$\Im\,\nabla\psi_0$}\}.
\end{equation}
The corresponding solution, that we call the {\it Legendre
transform potential}, can be written explicitly as,
\begin{equation}
\label{IVPToricGeodFormula}
\psi(s,x)=\psi_s(x):= (u_0 + s\dot{u}_0)^\star(x),
\quad (s,x)\in [0,T_\span^\cvx)\times\RR^n.
\end{equation}

The standard proofs that $\psi$ solves the HRMA
when $u=u(s,\,\cdot\,)$ solves the linear equation 
$\ddot u=0$, break down as soon as $\psi$ is not twice differentiable.
The underlying reason is that
the Legendre transform is not a bijection on non-convex functions.
Hence, as soon as $u_s$ ceases to be convex 
on the polytope $P$, one may 
not deduce---at least not using the classical 
reasoning---that $\psi$ solves the HRMA.  

Nevertheless, one may consider the formula (\ref{IVPToricGeodFormula}) as
defining a convex function (still denoted by $\psi$) on all of
$\RR_+\times\RR^n$. 
A natural question is whether $\psi$ continues to solve the HRMA
(\ref{HRMARayEq}) for $T>T_\span^\cvx$ in a weak sense. 

The main
result  of this article is that $\psi$ does not solve the HRMA
(\ref{HRMARayEq}) for any $T>T_\span^\cvx$. At the same time, we
prove that $\psi$ does  solve the equation wherever it
is differentiable, and in particular on a dense set whose
complement has zero Lebesgue measure. To state
our result, let
$$
\Delta(\psi):=
\{\,(s,x)\,:\, \psi \h{\rm\  is finite and differentiable at $(s,x)$\,}\}
\subset\RR_+\times\RR^n,
$$
denote the regular locus of $\psi$. Similarly, we denote the regular
locus of $\psi_s$ by $\Delta(\psi_s)\subset\RR^n$.
Let
$$
\Sigma_\sing:=\RR_+\times\RR^n\;\sm \Delta(\psi),
$$
denote the singular locus of $\psi$,
and set
$$
\Sigma_\sing(T):=\Sigma_\sing\cap [0,T]\times\RR^n.
$$
Since $\psi$ is everywhere finite the singular locus 
of $\psi$ has 
Lebesgue measure zero, and the regular locus is dense 
in $\RR_+\times\RR^n$ 

\begin{maintheo}
\label{FirstMainThm} 
Let $\psi$ be defined by (\ref{IVPToricGeodFormula}) for all $(s,x)\in\RR_+\times\RR^n$. 
Then

\noindent
(i) $\psi$ solves the HRMA (\ref{HRMARayEq}) on the regular locus. Namely,
$$
\MAop\psi=0 \h{ on } \Delta(\psi).
$$
In addition, $\,[0,T_\span^\cvx)\times\RR^n\subset \Delta(\psi)$.

\noindent
(ii) Whenever $T>T_\span^\cvx$, $\psi$ fails to solve the HRMA (\ref{HRMARayEq}).
In particular, the Monge-Amp\`ere measure of $\psi$ charges the
set $\Sigma_\sing(T)$ with positive mass and we have,
$$
\int_{[0,T]\times\RR^n} \MAop \psi
=
\int_{\Sigma_\sing(T)} \MAop \psi>0.
$$
\end{maintheo}

Let us outline the main steps in the proof.
In addition, a concrete overview of the proof is given
in \S\ref{OneDSection} for the case $n=1$.


In order for a convex function to be a weak solution of the HRMA the image of its
subdifferential must be a set of Lebesgue measure zero.
Our goal is therefore to obtain some description
of the subdifferential mapping $\del\psi$. First,  we study the regularity of 
the restriction $\psi_s$ of $\psi$ to each time slice.
Let
\begin{equation}
\label{AsFirstDefEq}
A_s:=\{y\in P\,:\, u_s(y)\ne\usd(y)\}\subset P.
\end{equation}

\begin{mainprop}
\label
{FirstMainProp}
For each $s>0$, the function
$$
\psi_s(x):= (u_0 + s\dot{u}_0)^\star(x),\quad x\in \RR^n,
$$
is a continuous strictly convex function on $\RR^n$. It
is Lipschitz continuous but not everywhere differentiable.
The singular 
locus of $\psi_s$ is given by
$$
\nabla \usd (A_s\sm\del P).
$$
\end{mainprop}

Geometrically, Proposition \ref{FirstMainProp} implies that $\psi$ can
be regarded as an infinite ray in the space of Lipschitz continuous
open-orbit \K potentials. 
Its proof is completed at the end of \S\ref{PartialRegularityLegendrePotentialSection},
and is divided into several steps. We also
remark that the strict convexity is not 
directly used for the proof of Theorem \ref{FirstMainThm}, however
several of the ingredients in its derivation are.

First, in order to prove that $\psi_s$ is strictly convex 
we study the regularity of its dual $\usd$.
We prove interior $C^1$ regularity for $\usd$ 
(Lemma \ref{LegendreLemma}). 
Then we prove that $\usd$ is 
moreover essentially smooth (Lemma \ref{EssentialSmoothnessLemma}), 
i.e., its gradient
blows-up on $\del P$. 
These results, together with a classical duality result,
then imply the strict convexity of $\psi_s$.
Here one uses the fact that $\psi_s$ is everywhere finite, and so
its subdifferential has domain equal to $\RR^n$, and also
$\Im\,\del\usd=\RR^n$. These facts 
can be seen directly for our explicit $\psi_s$, but
we also include a different proof (Lemma \ref{SurjectivityRnLemma})
using a homotopy argument that may have its own interest, 
and also shows as a by-product that $\usd$ and $\nabla\usd$
are continuous in $s$, a fact that is used later.

Second, to show that $\psi_s$ is not differentiable we show
that $\usd$ is not strictly convex. To that end,
for each $s>T_\span^\cvx$, we consider the set $A_s\subset P$ defined
in (\ref{AsFirstDefEq}). We show that this
set can be partitioned into convex sets 
$$
\overline{A_s}\sm\del P=\bigcup_{v\in A_s}Q(s,v)
$$
along which $\nabla\usd$ is constant 
(\S\ref{InvertibilitySection} and \S\ref{MAMeasureLegendrePotentialSection}).
Each such convex set 
$Q(s,y)$ (see (\ref{QMaxSYDefEq})) contains at least a line, and moreover if we
let $x=\nabla\usd(y)$ then $\del_x\psi(s,x)=Q(s,y)$
(Lemma \ref{YiPartialGammaPsiLemma}), 
proving that $\psi_s$ is not differentiable at $x$. 
We then obtain the exact description of the regular 
locus of $\psi_s$ in terms of $u_s$ and $A_s$ 
(Lemma \ref{FirstConeRegularityPsiLemma}), and this
concludes the proof of Proposition \ref{FirstMainProp}.

As just described, our results on the functions $\usd$ 
and their gradient maps already give a precise description
of the singularities of each $\psi_s$.
Next, we prove a partial $C^1$ regularity result 
for $\psi$ (Lemma \ref{PartialConeRegularityLemma})
that, as a corollary, gives a precise description of the 
singularities of $\psi$.
Recall that $\Delta(f)$ denotes the set on which $f$
is finite and differentiable.

\bigs\bigs

\begin{mainprop}
\label{SecondMainProp}
Assume that $x\in\Delta(\psi_s)$. Then $(s,x)\in\Delta(\psi)$.
I.e., the singular locus of $\psi$ is
the (indexed) union over $s$ of the singular loci of the 
functions $\psi_s$. In terms of the regular loci,
$$
\Delta(\psi)=\bigcup_{s\ge T_\span^\cvx}\{s\}\times\Delta(\psi_s).
$$
\end{mainprop}

In other words, wherever $\psi$ is
differentiable in $x$ it is also differentiable in $s$.
The results described so far then imply an
alternative description for the regular locus $\Delta(\psi)$
in terms of the maps $\del u_s,\, s\in\RR_+$ 
(Proposition \ref{SigmaRegCharProp}).
This 
allows us to show that $\psi$ solves
the HRMA on the regular locus. In fact, the image of
the (total) sub-differential of $\psi$ evaluated on the regular
locus is just the graph of $-\dot u_0$ over $P\sm\del P$, and this,
as a set in $\RR^{n+1}$, has Lebesgue measure zero
(Proposition \ref{MAPsiRegPartProp}).

To conclude the proof it thus remains to analyze the
\MA measure $\MAop\psi$ on the singular locus.
First, using some elementary facts regarding partial
subdifferentials of convex functions of several variables
we give a description of the $x$ partial subdifferential
of $\psi$ and of the set of reachable $x$ partial subgradients
of $\psi$ in terms of $u_s$ and the sets $Q(s,v)$ in the 
partition of $A_s$ (Lemma \ref{YiPartialGammaPsiLemma}).
Then, by using the partial regularity of $\psi$ we obtain some 
lower and upper bounds on $\del\psi(\{s\}\times\RR^n)$
(Lemma \ref{YyPartialGammaPsiLemma}). 

\begin{mainlem}
\label
{SecondMainLemma}

Let $y\in A_s\sm\del P$, 
and let $x:=\nabla\usdoublestar(y)$.
Then
\begin{equation*}
\h{\rm co}\,\big\{(-\dot u_0(v),v)\,:\, v\in \gamma_x\psi(s,x)\big\}
\subset
\del\psi(s,x)
\subset
\h{\rm co}\,\big\{(-\dot u_0(v),v)\,:\, v\in\del_x\psi(s,x)\big\}.
\end{equation*}
\end{mainlem}

 These bounds
are obtained in terms of certain convex sets projecting
onto the pieces $Q(s,v)$ of the partition of $A_s$. 
Using these bounds, along with monotonicity and continuity
of the family of the one-parameter family of sets $A_s$
(Lemmas \ref{AinfinityLemma} and \ref{AsContinuousMappingLemma}),
 we show that the subdifferential of $\psi$ 
in fact ``fills-in" a portion of the region lying between the graph 
of $-\dot u_0$ and the graph of minus the
convexification of $\dot u_0$ (see Figures 3 and 4), hence the mass of
$\MAop\psi$ necessarily becomes positive for any $T>T_\span^\cvx$
(Proposition \ref{MAMassPsiProp} and Lemma \ref{MASingMassAinfinityLemma}),
completing the proof. 

It is worth pointing out that the proof also shows that
the Monge-Amp\`ere mass of $\psi$ has an a priori upper bound 
depending only on the Cauchy data (Lemma \ref{MASingMassAinfinityLemma}).
Let $\h{\rm epi}\, f$ denote the epigraph of $f$ 
(see \S\S\ref{ConvexAnalytisSubsection}).

\begin{mainprop}
\label{ThirdMainProp}
Let $T>0$. One has,
$$
\int_{[0,T]\times\RR^n}\h{\rm MA}\,\psi
\le\h{\rm Vol}\, 
\big(\,
 \epi\,(-\dot u_0)\, \sm \epi\,(-(\dot u_0)^{\star\star})
\,\big),
$$
where the right hand side denotes the volume in $\RR^{n+1}$ of
the set of points lying below the graph of $\dot u_0$ and above the graph
of its convexification, over $P$.
\end{mainprop}

The graph of $\psi$ defines a hypersurface in $\RR^{n+2}$
that is flat over $[0,T_\span^\cvx)\times\RR^n$. Proposition
\ref{ThirdMainProp} implies that the non-compact surface in $\RR^{n+2}$ 
defined by the graph of $\psi$, while not flat, 
has finite total Gauss curvature (see 
\S\S\ref{OtherCauchySubsection} and  \S\ref{OneDSection}).
 
Note that the right hand side is zero if and only if $\dot u_0$ is convex.
The inequality is also sharp when the right hand side is positive.
This can be seen by considering, for instance, explicit examples 
in $n=1$ (see \S\ref{OneDSection}), where equality is attained 
in the limit where $T$ tends to infinity.

It is interesting to give a geometric description of the singular
locus $\Sigma_\sing$ on which the \MA mass is concentrated.
By Propositions \ref{FirstMainProp} and \ref{SecondMainProp} we have
an explicit description of the singular locus,
\begin{equation}
\label{SigmaSingDescriptionEq}
\Sigma_\sing
\,=
\bigcup_{s\ge T_\span^\cvx,}\{s\}\times\nabla \usd(A_s\sm\del P).
\end{equation}
Now, by the continuity of $\nabla\usd$ in $s$ and
the set-valued continuity of the sets $A_s$,
it follows that $\Sigma_\sing$ is a countable union of $C^0$ hypersurfaces in $\RR^{n+1}$.
Moreover, by general results of Alberti \cite{A} in fact $\Sigma_\sing$ must
be a countable union of locally Lipschitz continuous hypersurfaces 
in $\RR^{n+1}$. Further, for reasonable Cauchy data, e.g., such that 
$A_s$ has finitely many components for all $s>0$, it follows from
(\ref{SigmaSingDescriptionEq}) that $\Sigma_\sing$
will also be composed of finitely many hypersurfaces.
A visualization of the 
singular set and the corresponding `corner set' of the graph of $\psi$
is given in \S\ref{OneDSection} (see Figures 2 and 3).

\subsection{Consequences for the quantum lifespan}
\label{QuantumConsequencesSubsection}

We now tie Theorem \ref{FirstMainThm} together with the main
result of \cite{RZ2}. There we constructed a candidate solution 
of the IVP for the HCMA on a general projective \kahler manifold
$(M,\o)$ by a quantization procedure, inspired by
the formal analytic continuation argument of Semmes and 
Donaldson \cite{D1,S2} and by the Phong-Sturm construction
of geodesic segments \cite{PS1} and test configurations geodesic 
rays \cite{PS2} by finite-dimensional approximations.  
We first quantized the Hamiltonian flow
determined by the initial velocity $\dot{\varphi}_0$
and metric $\omega_{\varphi_0}:=\o+\i\ddbar\vp_0$
as a semi-classical one-parameter subgroup of unitary Toeplitz
complex Fourier integral operators

$$
U_N(t):=\Pi_N e^{\i tN \Pi_N\dot\vp_0\Pi_N}\Pi_N
$$
\smallskip
\noindent on the spaces  $H^0(M,L^N)$ of holomorphic sections of the
quantizing line bundle $L \to M$. We then analytically continued
the unitary group to an imaginary time semi-group
\begin{equation}
\label{UKDEF} U_N(\i s)  : H^0(M, L^N) \to H^0(M, L^N).
\end{equation}
We denote by $U_N(-\i s)(z,w)$ the Schwartz kernel of this
operator with respect to the volume form $(N\omega_{\vp_0})^n$.
The potentials \label{QuantumAnalyticPotentialDef}
\begin{equation}
\label{NAnalyticContPotential} \vp_N(s,z):= \frac{1}{N} \log
U_N(-\i s,z,z)
\end{equation}
are then readily seen to be subsolutions of the HCMA

\begin{equation}
\label{HCMARayEq}
\left\{
\begin{array}{rrl}
(\pi_2^\star\omega + \i\ddbar \vp)^{n+1}
\!\!\!& = & \!\!\!
0 \quad\quad\;\,\mskip2mu \mbox{on} \; S_{T} \times M,
\cr\cr
\vp(0,s,\,\cdot\,)
\!\!\!& = &\!\!\!
\vp_0(\,\cdot\,)
\;\; \mbox{on} \; \{0\}\times\RR \times M,
\cr\cr
\displaystyle
\frac{\partial\vp}{\partial s}(0,s,\,\cdot\,)
\!\!\!& = &\!\!\!
\dot\vp_0(\,\cdot\,) \;\; \mbox{on} \; \{0\}\times\RR \times M,
\end{array} \right.
\end{equation}
where $\pi_2:S_T\times M\ra M$ is the projection (recall $S_T=[0,T]\times\RR$).
We then defined the {\it quantum analytic
continuation potential} $\vp_\infty$ by
\begin{equation} 
\label{QPotentialEq}  
\vp_\infty(s,z) := \lim_{l \to \infty} (\sup_{N \geq
l}\vp_N)_\uscreg(s,z),
\end{equation}
where 
$u_\uscreg(z_0)
:= 
\lim_{\epsilon \to 0} \sup_{|z - z_0|<\epsilon} u(z)
$ 
denotes the upper semi-continuous regularization of $u$.
The {\it quantum lifespan} was then defined as 
$$
T_\span^Q:=\sup\,\{\,T\,:\, \vp_\infty \h{ solves (\ref{HCMARayEq}) on } S_T\times M\,\}.
$$

In the toric case, the Cauchy data $\o_{\vp_0},\,\dot\vp_0$ gives rise 
to Cauchy data $\psi_0,\,\dot\psi_0$ for (\ref{HRMARayEq}), 
where $\o_{\vp_0}$ equals $\i\ddbar\psi_0$ on the open orbit,
and $\dot\psi_0$ simply denotes the restriction of $\dot\vp_0$ to the open orbit.
We proved in  \cite{RZ2}  
that the quantum
analytic continuation potential coincides, on all of $\RR_+\times\RR^n$,
with the Legendre transform
potential:
$$
\psi_0+\vp_\infty(s,\,\cdot\,)=\psi_s
$$ 
(by the symmetry, $\vp_\infty$ can be regarded as a function on $\RR_+\times\RR^n$).
Theorem \ref{FirstMainThm} therefore implies:

\begin{maincor} 
\label{ACPOT} 
On a toric or Abelian variety, $T_\span^Q=T_\span^\cvx$.
In other words, the quantum analytic continuation
potential $\vp_\infty$ given by
(\ref{QPotentialEq}) is a subsolution of the HCMA (\ref{HCMARayEq}) 
that solves the Cauchy problem 
until $T_\span^\cvx$,  and cease to solve it, even as a
weak solution, after that time.
\end{maincor}

As discussed in \cite{RZ2}, \S1 and \S\S3.1, 
the question whether there exist alternative
solutions of HRMA and HCMA is taken up in sequels to
this article. In particular, in \cite{RZ3} we 
give a characterization of the smooth lifespan $T_\span^\infty$ 
of the HCMA (\ref{HCMARayEq}) from which we conclude
that there exists no smooth solution to the HRMA
(\ref{HRMARayEq}) beyond the convex lifespan, 
namely $T_\span^\cvx=T_\span^\infty$.

\subsection{Examples on $S^2$}
\label
{ExampleStwoSubsection}

We illustrate with an example the behavior of the family of \K metrics (viewed
as a path in $\H_\o$)
associated with the quantum analytic continuation (or Legendre transform)
potential.

Let $(\PP^1,\o)$ denote the Riemann sphere equipped with a \K form $\o$.
Consider an $S^1$-invariant \K metric $\o_{\vp_0}$ on $\PP^1$
that equals $\i\ddbar\psi_0$ away from the poles.
Let $\dot{\vp}_0$ denote a given $S^1$ invariant initial velocity.
The geodesic equation
\begin{equation}
\label{GeodesicTwoDimEq}
\frac{\del^2\psi}{\del s^2}
=
\Big(\frac{\del^2\psi}{\del x^2}\Big)^{-1}
\Big(\frac{\del^2\psi}{\del x\del s}
\Big)^2,
\end{equation}
can be interpreted as the HRMA $\det\nabla^2\psi=0$.
Denote 
$$
\psi'(s,x):=\frac{\del\psi(s,x)}{\del x},\quad
\dot\psi(s,x):=\frac{\del\psi(s,x)}{\del s}.
$$
By letting $y:=\psi'(s,x)$, and
$u(\psi'(s,x)):=x\psi'(s,x)-\psi(s,x)$, a computation
shows that equation (\ref{GeodesicTwoDimEq})
becomes $\ddot u(s,y)=0$, 
solved by $u(s,y):=u_0(y)+s\dot u_0(y)$, with
$\dot u_0(y)=-\dot\psi_0((\psi_0')^{-1}(y))$.

Note that $y=y_s,\th$ are the action-angle variables
for $\omega_s=\i\ddbar\psi(s,x)$, i.e.
$\omega_s = dy_s \wedge d\theta$
and $y_s$ is the moment map of the $S^1$ action with respect to
$\omega_s$.
Let $z=e^{x+\i\th}$ denote the holomorphic
coordinate away from the poles. The metric $g_s$ at time $s$ is
then given by
\begin{equation}
\label{GsZerothEq}
g_s=\psi''(s,x)\big(dx^2+d\th^2),
\end{equation}
expressed in action-angle variables as
\begin{equation}
\label{GsFirstEq}
g_s
=
u(s)_{yy} dy^2 + \frac{1}{u(s)_{yy}} d \theta^2,
\end{equation}
where $u(s)_{yy}$ denotes the second derivative of the symplectic
potential $u$ at time $s$ with respect to the action variable
$y=y_s$. Also, let $r_s$ denote the geodesic distance function from
the north pole (the fixed point of the $S^1$ action on which $y_s$
takes its maximum). Then $r_s$ is a function only of $y_s$.
Hence the change of variables from $y_s$ to $r_s$ does not add any
term containing $d \theta$, and
\begin{equation}
\label{GsSecondEq}
g_s
=
dr_s^2 +  \frac{1}{u(s)_{yy}} d \theta^2.
\end{equation}

Formulas (\ref{GsFirstEq})--(\ref{GsSecondEq}) are valid only
so long as $u(s,\,\cdot\,)=u_0+s\dot u_0$ remains convex,
and they show that in that regime $g_s$ is a smooth metric.
The Legendre transform potential $\psi$ defined by
(\ref{IVPToricGeodFormula}) provides an extension
of the path of metrics $\{g_s\}_{s\in[0,T_\span^\cvx)}$ given
by (\ref{GsFirstEq}) to $s\in[T_\span^\cvx,\infty)$, given by
\begin{equation}
\label{GsThirdEq}
g_s={u_s^\star}''(s,x)\big(dx^2+d\th^2).
\end{equation}
Equation (\ref{IVPToricGeodFormula}) also means that $\psi(s,\,\cdot\,)$ is the
Legendre dual of the convexification of $u(s,\,\cdot\,)$.
The convexified symplectic potential $\usdoublestar$ is $C^1$ but
develops straight segments on intervals of non-convexity (see Figure 1).
At these, $u^{\star\star}(s)_{yy} = 0$. The radial
`height' of the inverse image of this straight segment under $y_s$
is zero but the lattice circle ($S^1$ orbit) is of infinite
radius. Thus, the metric develops an $S^1$-invariant
delta-function singularity at the corresponding value of $r_s$.
Theorem \ref{FirstMainThm} is the statement that the path
of metrics $g_s$ ceases to be a geodesic precisely at
$T_\span^\cvx$, when the singularities appear. However,
it does satisfy the geodesic equation on a dense set,
whose time slice is the complement of a discrete set of singular $S^1$ orbits.

\bigskip\bigskip\bigskip

\hglue1.1in\includegraphics[scale=0.5]{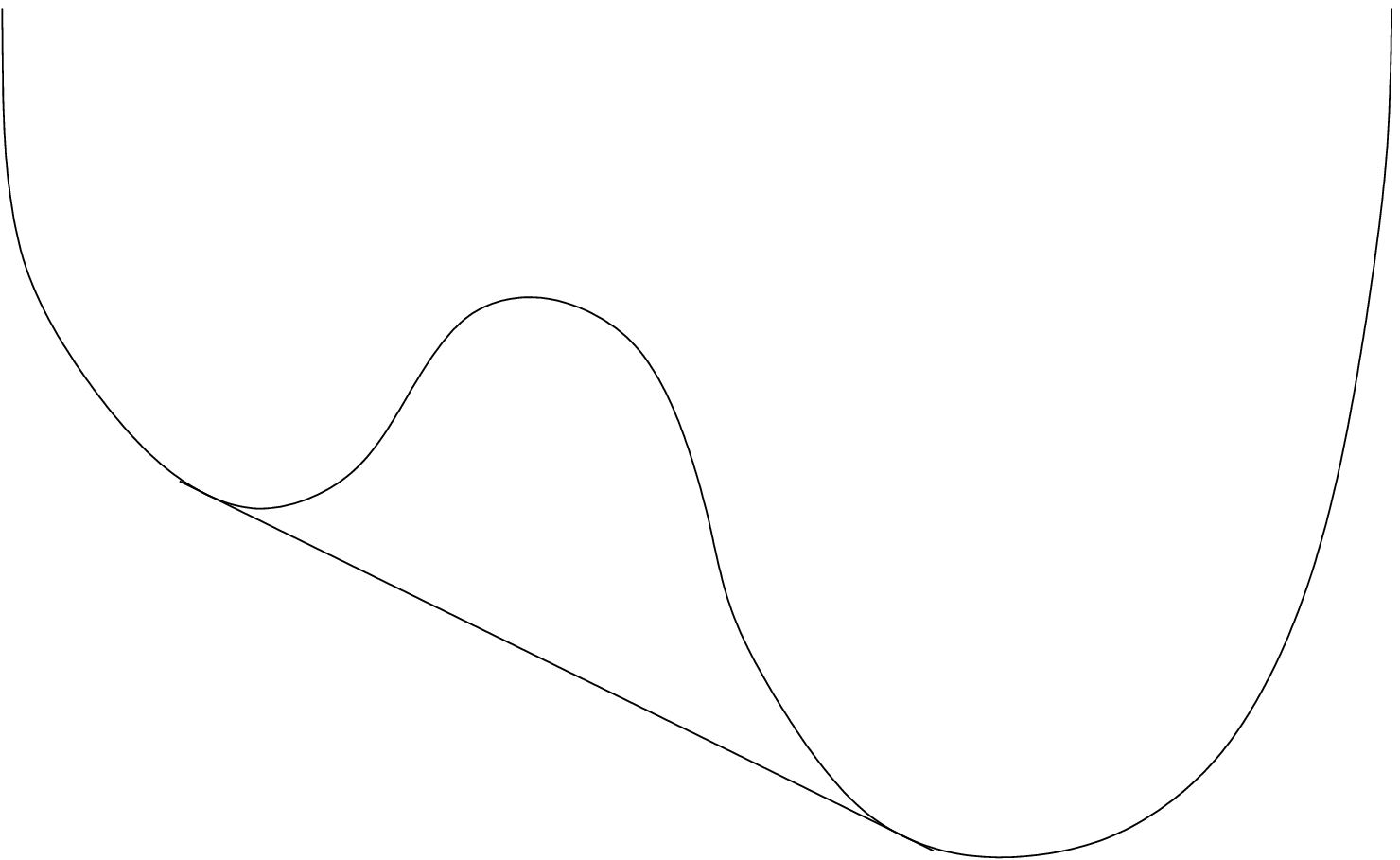}


\bigskip
\centerline{
{\sc Figure 1.} The graphs of $ u_s$ and $\usdoublestar$
over $P$ for $s>T_\span^\cvx$.
}

\bigskip\bigskip

\subsection{Asymptotic behavior of geodesic and subgeodesic rays}
\label{AsymptoticGeodesicRaysSubsection}

It has been conjectured by Donaldson that smooth geodesic
rays in the space of \K metrics should play an important 
role in questions regarding geometric `stability' and 
existence of canonical metrics \cite{D1}. 
An important question is therefore precisely which directions
yield smooth geodesic rays.
As mentioned in \S\S\ref{QuantumConsequencesSubsection}, 
we prove in the sequel that for the HRMA 
(\ref{HRMARayEq}) one has $T_\span^\infty=T_\span^\cvx$.
Therefore the directions of smooth toric geodesic rays
are precisely those with infinite convex lifespan.
Further, as a corollary of the results here and in the sequel 
we describe the limiting behavior of the rays obtained by solving 
the IVP using the quantum method of \cite{RZ2}, or, equivalently,
by means of the Legendre transform method. This holds whether they be 
smooth geodesics or only `subgeodesics', by which we mean 
subsolutions of the HCMA.

\subsection{Other Cauchy problems for the HRMA}
\label
{OtherCauchySubsection}

The IVP for geodesics in the space of toric metrics on toric varieties 
gives rise to a specific Cauchy problem for the HRMA. A natural question
is to what extent do the techniques used in this article  
generalize to other Cauchy problems for the HRMA. We briefly discuss
several possible problems and generalizations, some of which we hope
to discuss in more detail elsewhere.

Aside from the fact that
the Cauchy data and the Cauchy hypersurface are smooth,
and that we require all solutions to be convex,
the main distinctive features of the
Cauchy problem (\ref{HRMARayEq}) are:

\medskip

(i)
the Cauchy hypersurface is an $\RR^n$-slice of the total space $[0,T]\times\RR^n$,

\smallskip
(ii)
the initial convex function has linear growth at infinity, and the\hfill\break
\hglue33.3pt initial velocity is uniformly bounded.

\medskip

An example of a Cauchy problem where (i), but not (ii), holds, is the IVP
for geodesics in the space of toric metrics on Abelian varieties.
In that setting the convex function $\psi_0$ has quadratic growth
at infinity. On the other hand, the gradient map $\nabla\psi_0$
now maps to a torus instead of a polytope. It follows that the
symplectic potential $\psi_0^\star$ can be considered as 
a convex function on $\RR^n$ with quadratic growth and no singularities.
This fact eliminates the need for the analysis near the boundary of $P$
that was necessary in this article. Although we do not go into 
the details here, using the methods of this series one may prove
an analogue of our main results for this class of manifolds.

Note that the requirement in (ii) that the initial 
velocity be bounded, while the initial convex 
function divergent at infinity, guarantees that 
the gradient image of $\psi_s$ remains constant, and this
still held true in the Abelian case.
Removing this assumption from (ii) would require dealing
with a one-parameter family of gradient images $P_s$,
and would certainly complicate some of the analysis.

Next, one may allow the Cauchy hypersurface to be more general,
for example a smooth affine hypersurface in an affine manifold.
Geometrically, this setting arises when one considers
the IVP for geodesics in spaces of invariant metrics
in various classes of manifolds with large symmetry 
groups (for examples see, e.g., \cite{D2}, \S4).
It would be interesting to obtain a formula 
for $T_\span^\infty$ analogous to
the one for the convex lifespan (\ref{TspanCvxEq}).
We note that in the affine setting Foote \cite{F1,F2} has
given a sufficient condition on the Cauchy data and hypersurface
to have $T_\span^\infty>0$.  

Although the affine situation is certainly more complicated
than the Euclidean one it seems plausible to us that one could 
generalize at least some of the techniques of the present series 
to the affine setting.

Finally, we mention that the Cauchy problem for the HRMA is also
classically used to construct smooth Gauss flat hypersurfaces in $\RR^{n+2}$
(this is carried out in $\RR^3$ in \cite{U}).
It seems interesting to investigate what kind of singular hypersurfaces
are obtained from the Legendre transform method when extended
beyond the convex lifespan. A consequence of
the a priori bound of Proposition \ref{ThirdMainProp} is
that the resulting hypersurface has bounded total Gauss
curvature, a fact that seems rather surprising.
In \S\ref{OneDSection} we
briefly touch upon this point of view by proving Theorem
\ref{FirstMainThm} in the case $n=1$, and also illustrate
explicitly the finite curvature phenomenon.

\bigskip
\section{Background}
\label{SectionBackground}

\bigs

In this section, we recall some basic definitions relating to
convex analysis, the real \MA operator, and the reduction of HCMA on 
toric \kahler manifolds to the HRMA. 
For additional necessary notation and definitions that are used 
throughout we refer to \cite{RZ2}, \S\S4.3.

\subsection{Convex analysis}
\label{ConvexAnalytisSubsection}

For general background on Legendre duality and convexity we refer
the reader to \cite{HL1,HL2,Ro}, whose notation and terminology
we largely adhere to. Denote by
$$
\h{\rm co}\, A
$$ the convex hull of the set $A$. Given a
function $f$ defined on a set $P\subset\RR^n$ let
$$
\epi\,f:=\{(x,r)\,:\, r\ge f(x), \, x\in P\}
$$
denote its epigraph. A vector
$v\in(\RR^n)^\star$ is said to be a subgradient of a function $f$
at a point $x$ if $f(z)\ge f(x)+\langle v,z-x\rangle$ for all $z$.
The set of all subgradients of $f$ at $x$ is called the
subdifferential of $f$ at $x$, denoted $\partial f(x)$.

The Legendre-Fenchel conjugate of a continuous function 
$f=f(x)$ on $\RR^n$ is defined by \cite{Fe}
$$
f^\star(y):=\sup_{x\in\RR^n}\big(\langle x,y \rangle - f(x)\big).
$$
For simplicity,  we will refer to $f^\star$ sometimes as the Legendre dual, or just dual,
of $f$.

\begin{definition}
\label{ConvexTermsDef}
A convex function $f$ is called proper if it is not identically $+\infty$ and is uniformly
bounded below. A proper convex function is called closed if it is lower semi-continuous.
The domain of $f$ is defined by $\dom(\partial f):=\{x\,:\, \partial f(x)\ne\emptyset\}$.
\hfill\break
(i)
A function on a convex set $C$ is called strictly convex on $C$ if
for every $\lambda\in(0,1)$ and all distinct points $x_1,x_2\in C$
holds
$f((1-\lambda)x_1+\lambda x_2)<(1-\lambda)f(x_1)+\lambda f(x_2)$.
\hfill\break
(ii)
A proper convex function is called essentially strictly convex if it is strictly convex
on every convex subset of $\dom(\partial f)$.
\hfill\break
(iii)
A proper convex function is called essentially smooth if $C:=\h{\rm int}(\dom(\partial f))\ne\emptyset$,
if $f$ is differentiable on $C$ and if $\lim_{i\ra\infty}|\nabla f(x_i)|=+\infty$
whenever $\{x_i\}_{i\ge1}$ is a sequence in $C$ converging to $x\in\partial C$.
\end{definition}

\subsection{The real Monge-Amp\`ere operator}
\label{MAOpSubsection}

We recall the definition of the Monge-Amp\`ere operator
and its basic characterization, due to Alexandrov and Rauch-Taylor.
Let $M(\RR^{n+1})$ denote the space of differential forms of degree $n+1$ on $\RR^{n+1}$
whose coefficients are Borel measures (i.e., currents of degree $n+1$ and order 0).

\begin{prop} {\rm (See \cite{RT}, Proposition 3.1)}
Define by 
$$
\h{\rm MA} f:=
d\frac{\partial f}{\partial x^1}\w \cdots\w d\frac{\partial f}{\partial x^{n+1}},
$$
an operator $\h{\rm MA}: C^2(\RR^{n+1})\ra M(\RR^{n+1})$.
Then $\h{\rm MA}$ has a unique extension to a continuous operator on the cone
of
convex functions.
\end{prop}

An alternative, geometric, definition is due to Alexandrov, and uses the notion
of a subdifferential of a convex function.

\begin{prop}{\rm (See \cite{RT}, Section 2)}
\label{RTProp}
For any convex function $f$, the measure $\calM\calA\, f$, defined
by
$$
(\calM\calA\, f)(E):=\h{\rm Lebesgue measure of\ } \partial f(E),
$$
is a Borel measure.
\end{prop}

The following result of Rauch-Taylor links these two definitions and will be crucial below.

\begin{theo} {\rm (See \cite{RT}, Proposition 3.4)}
\label{RTThm}
For every convex function $f$ on $\RR^{n+1}$ one has the equality
of Borel measures $\h{\rm MA} f = \calM\calA\, f$.
In particular, the real Monge-Amp\`ere measure is zero if and only if
the image of the subdifferential map has Lebesgue measure zero in $\RR^{n+1}$.
\end{theo}

\bigskip

\subsection{The Cauchy problem for the symplectic potential}
\label{SectionCauchySympPot}

It is well-known that the Legendre transform linearizes the HRMA.
This fact also has a geometric interpretation that we now briefly review.

Let $(M, \omega)$ be a toric \kahler manifold of complex dimension
$n$ and let $\bfT=(S^1)^n$ denote the real torus of dimension $n$ which acts
on $(M, \omega)$ in a Hamiltonian fashion. We denote by
$\hcal(\bfT)$ the class of $\bfT$-invariant \kahler metrics in the
cohomology class of $\omega$. On the  open-orbit of 
$\bfT^{\CC}=(\CC^\star)^n$, a $\bfT$-invariant 
\kahler metric has a \K potential $\psi$ and we
also write  $\psi\in\calH(\bfT)$.

Since it is $\bfT$-invariant, the \K potential may be identified
with a smooth strictly convex function on $\RR^n$ in logarithmic
coordinates. Therefore its gradient $\nabla\psi$ is one-to-one
onto $P=\overline{\Im\nabla\psi}$ and one has the following explicit
expression for its Legendre dual (\cite{Ro}, or \cite{R}, p.
84--87),
\begin{equation}
\label{LegendreDualityEq} u(y)=\psi^\star(y)=\langle
y,(\nabla\psi)^{-1}(y)\rangle-\psi\circ(\nabla\psi)^{-1}(y),
\end{equation}
which is a smooth strictly convex function on $P$, satisfying
\begin{equation}
\label{LegendreDualityGradientEq} \nabla
u(y)=(\nabla\psi)^{-1}(y),
\end{equation}
and
\begin{equation}
\label {LegendreDualityHessianEq} \big(\nabla^2
u(y)\big)^{-1}=\nabla^2\psi((\nabla\psi)^{-1}(y)).
\end{equation}
Following Guillemin, the function $u$ is called the
symplectic potential of the metric $\i\ddbar\psi$. The space of all
symplectic potentials is denoted by $\calL\calH(\bfT)$. Put
\begin{equation}
\label{GuilleminFormulaEq} u_G:=\sum_{k=1}^d l_k\log l_k.
\end{equation}
A result of Guillemin \cite{G1} states that for any symplectic
potential $u$ the difference $u-u_G$ is a smooth function on $P$
(that is, up to the boundary). In other words,
\begin{equation}
\label{LHTDefSecondEq} \calL\calH(\bfT) =\{u\in C^\infty(P\sm\del
P)\,:\, u=u_G+F,\quad\h{\rm with } F\in C^\infty(P)\}.
\end{equation}

The Legendre transform is an isometry between $(\H(\bfT),g_{L^2})$
and $(\calL\H(\bfT),L^2(P))$. It transforms the Christoffel
symbols of $(\H(\bfT),g_{L^2})$ to zero and thus linearizes the
Monge-Amp\`ere equation to the equation $\ddot
u=0$ (for more on this see \cite{Gu,S1}, or \cite{RZ1}, \S3). The
differential of the Legendre transform acts as minus the identity,
that is if $\eta_s$ is a curve in $\H(\bfT)$ and if $
u_s:=\eta^\star_s $ are the corresponding symplectic potentials
then
\begin{equation}
\label{VariationPotentialEq}
\dot\eta_s =-\dot u_s\circ \nabla\eta_s
\end{equation}
(see, e.g., \cite{R}, p. 85).
Therefore the IVP on $(\H(\bfT),g_{L^2})$ is transformed to the following initial value
problem for geodesics in the space of symplectic potentials:
\begin{equation}
\label{SympPotentialIVPEq}
\ddot u=0, \quad, u_0=\psi^\star_0,\quad \dot u_0=-\dot\psi_0\circ(\nabla\psi_0)^{-1}.
\end{equation}

\bigs
\section{Legendre continuation of flat surfaces in $\RR^3$}
\label{OneDSection}
\bigs

Our purpose in this section is to explain the proof of Theorem 
\ref{FirstMainThm} in the case where $n=1$, i.e., for
the Cauchy problem for the 2-dimensional HRMA
of \S\S\ref{ExampleStwoSubsection}.
This special setting is simpler to visualize explicitly
and helps motivate and better capture some of the
constructions carried out in the proof of Theorem
\ref{FirstMainThm}, and also, part of the analysis
simplifies. Many of the assertions in this section
are stated without rigorous justification, and their proofs
can be found (for all $n$) in later sections.

The graph of $\psi(s,x)$
$$
\{(s,x,\psi(s,x)\}\subset\RR^3
$$
over $[0,T_\span^\cvx]\times\RR$ is a flat surface 
since its second fundamental form is proportional
to the Hessian of $\psi$. We are interested in
what happens to this surface when extended beyond
$T_\span^\cvx$: does $\psi$ still define a flat
surface, in a weak sense, after $T_\span^\cvx$?

Let $\oFS$ denote the Fubini-Study form of constant 
Ricci curvature $1$ on the Riemann sphere, 
given locally by
$$
\oFS=\frac\i{\pi}\frac{dz\w\dbz}{(1+|z|^2)^2}.
$$
The associated open-orbit \K potential can be taken as 
\begin{equation}
\label{PsizeroStwoEq}
\psi_0(x)=\log(1+|z|^2)-\half\Re z=\log(1+e^{2x})-x.
\end{equation}
The corresponding moment polytope is $[-1,1]$, and 
the symplectic potential dual to $\psi_0$ can be computed
via the moment map $y(x)=\psi_0'(x)\in [-1,1]$, 
\begin{equation}
\label{UzeroStwoEq}
u_0(y)=(1+y)\log(1+y)+(1-y)\log(1-y), \quad y\in[-1,1].
\end{equation}
Let $\dot\vp_0\in C^\infty(S^2)$ be given and set
$$
\dot u_0=-\dot\vp_0((\psi_0')^{-1}(\,\cdot\,))
=-\dot\psi_0((\psi_0')^{-1}(\,\cdot\,))\in C^\infty([-1,1]).
$$

We start with the analysis of $\usd$. 
By (\ref{UzeroStwoEq}) we have 
$\lim_{y\ra\pm1}u_0'=\pm\infty$, and the
same holds for $u_s$ since $\dot u_0$ is 
bounded on $P$. Therefore, when
restricted to a neighborhood of $\{\pm1\}$ in $[-1,1]$
the tangent lines to the graph of $u_s$ lie
below the graph. Hence, 
$A_s\cap\{\pm1\}=\emptyset$, and the graphs of
$u_s$ and of $\usd$ differ only above 
$(-1+\epsilon,1-\epsilon)$ for some $\epsilon>0$.
Above every connected component of $A_s$ the graph
of $\usd$ is necessarily affine with slope precisely
equal to the derivative of $u_s$ at the end-points.
Hence $\usd$ is continuously differentiable.

Let $s>T_\span^\cvx$ so $A_s\ne\emptyset$.
Let us assume for simplicity that $A_s$ is connected,
with $A_s=(a_s,b_s)\subset[-1,1]$, and set $x_s:=\nabla\usd\big((a_s,b_s)\big)$.
Then the graph of $\psi_s$ has a corner 
of angle $\alpha_s=\tan^{-1}a_s-\tan^{-1}b_s$ 
at $(x_s,\psi_s(x_s))$, and is smooth elsewhere.
The set of reachable subgradients of $\psi_s$ at $x_s$ 
is $\gamma\psi_s(x_s)=\{a_s,b_s\}$. It follows that
$\gamma\psi(s,x_s)$, the set of reachable sugradients
of $\psi$ at $(s,x_s)$, contains 
$\{(-\dot u_0(a_s),a_s),(-\dot u_0(b_s),b_s)\}$. Thus,
by convexity, the total subdifferential of $\psi$
at $(s,x_s)$ contains the line connecting these two points:
\begin{equation}
\label{LineIntervalEq}
\h{co}\,\{(-\dot u_0(a_s),a_s),(-\dot u_0(b_s),b_s)\}
\subset\del\psi(s,x_s).
\end{equation}

On the other hand, the graph of $-\dot u_0$ 
cannot be linear even locally on $(a_s,b_s)$ since
that would imply that $u_0+s\dot u_0$ were
strictly convex on some part of $(a_s,b_s)$. It follows
that the affine segment in (\ref{LineIntervalEq})
intersects the graph of $-\dot u_0$ over $(a_s,b_s)$
only in a discrete set of points. Since $A_s$
is monotonically increasing in $s$ in a continuous 
fashion, we conclude that the image of $\del\psi$ 
fills-in a region in $\RR\times[-1,1]$ 
above the graph of $-\dot u_0$, and consequently that $\psi$
is not a weak solution of the 2-dimensional HRMA on
$[0,T]\times\RR$ for any $T>T_\span^\cvx$.
Moreover, the region filled out, as $T$ tends to infinity, 
is precisely the region above the graph of $-\dot u_0$ and below
the graph of $-(\dot u_0)^{\star\star}$.
We conclude that the non-compact surface defined by $\psi$ is flat precisely
on the complement of the corner set, but not globally flat.
However, its total curvature in the sense of Alexandrov while positive, is finite, 
and is concentrated on the codimension one corner set.

\vglue1.5cm

\hglue2.1cm\includegraphics[scale=0.5,angle=-90]{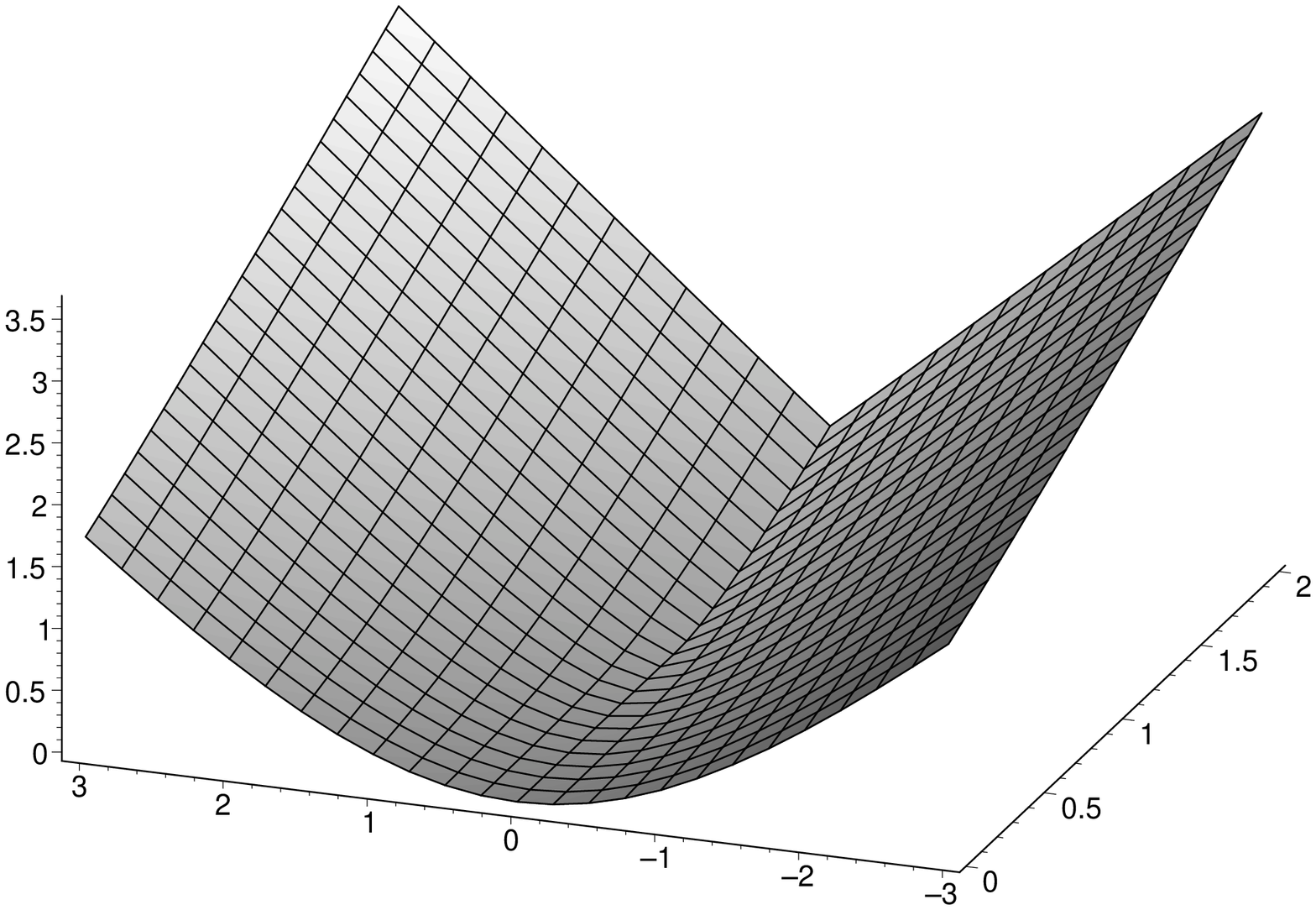}


\bigskip\bigskip
\noindent
{\sc Figure 2.} The graph of $\psi$ for $(s,x)\in[0,2]\times[-3,3]$
and Cauchy data $(\psi_0,\dot\psi_0)$ given by 
(\ref{PsizeroStwoEq}) and (\ref{PsiDotZeroEq}).

\vglue-5.5cm\hglue1.3cm$\psi$
\vglue5.5cm

\vglue-2.7cm\hglue5.3cm$x$
\vglue2.7cm

\vglue-4.3cm\hglue11.3cm$s$
\vglue4.3cm
\vglue.5cm

\begin{exam}
{\rm
The graph of $\psi$ may be parametrized with respect to the moment
coordinate on the domain of its invertibility, namely, 
$$
\h{\rm graph of $\psi\,$}
=
\Big\{
\big(\,s,\,\nabla u_s(y),\,\langle y,\nabla u_s(y)\rangle-u_s(y)\,\big)
\,:\, s\in\RR_+,\; y\in P\sm(A_s\cup\del P)
\Big\}.
$$
To visualize the graph of $\psi$ in a specific example
let $\psi_0$ be given by (\ref{PsizeroStwoEq}), and set 
\begin{equation}
\label{PsiDotZeroEq}
\dot\vp_0(z)
=
\bigg(\frac{|z|^2-1}{|z|^2+1}\bigg)^2,\quad \h{or,} \quad
\dot\psi_0(x)
=
\bigg(\frac{e^{2x}-1}{e^{2x}+1}\bigg)^2.
\end{equation}
Hence, $\dot u_0(y)=-y^2$.

A portion of this graph is drawn in Figure 2.
In this example, $T_\span^\cvx(\psi_0,\dot\psi_0)=1$, and the graph of $\psi$ develops a kink
along the corner set above 
$\Sigma_\sing=\{x=0, \, s>1\}.$ 
For each $s>1$ the graph of $\psi_s$ has a kink at $x=0$, and so the
corner set of $\psi$ is precisely the union over $s$ of the corner sets of
$\psi_s$. Moreover, in this special example $u_0$ and $\dot u_0$ are even functions. 
The same is then true for $\psi_s$, and one may further show that $\psi$ is
differentiable in $s$. Let $A_s=(-a_s,a_s)$. Then 
$\del_x\psi(s,x_s)=[-a_s,a_s]$, and 
$$
\del\psi(s,x_s)=\{\dot\psi_s(x_s)\}\times\del_x\psi(s,x_s)=\{(a_s^2,y)\,:\, y\in[-a_s,a_s]\}.
$$

\vglue1.5cm

\hglue1.7cm
\includegraphics
[scale=0.5,angle=-90]
{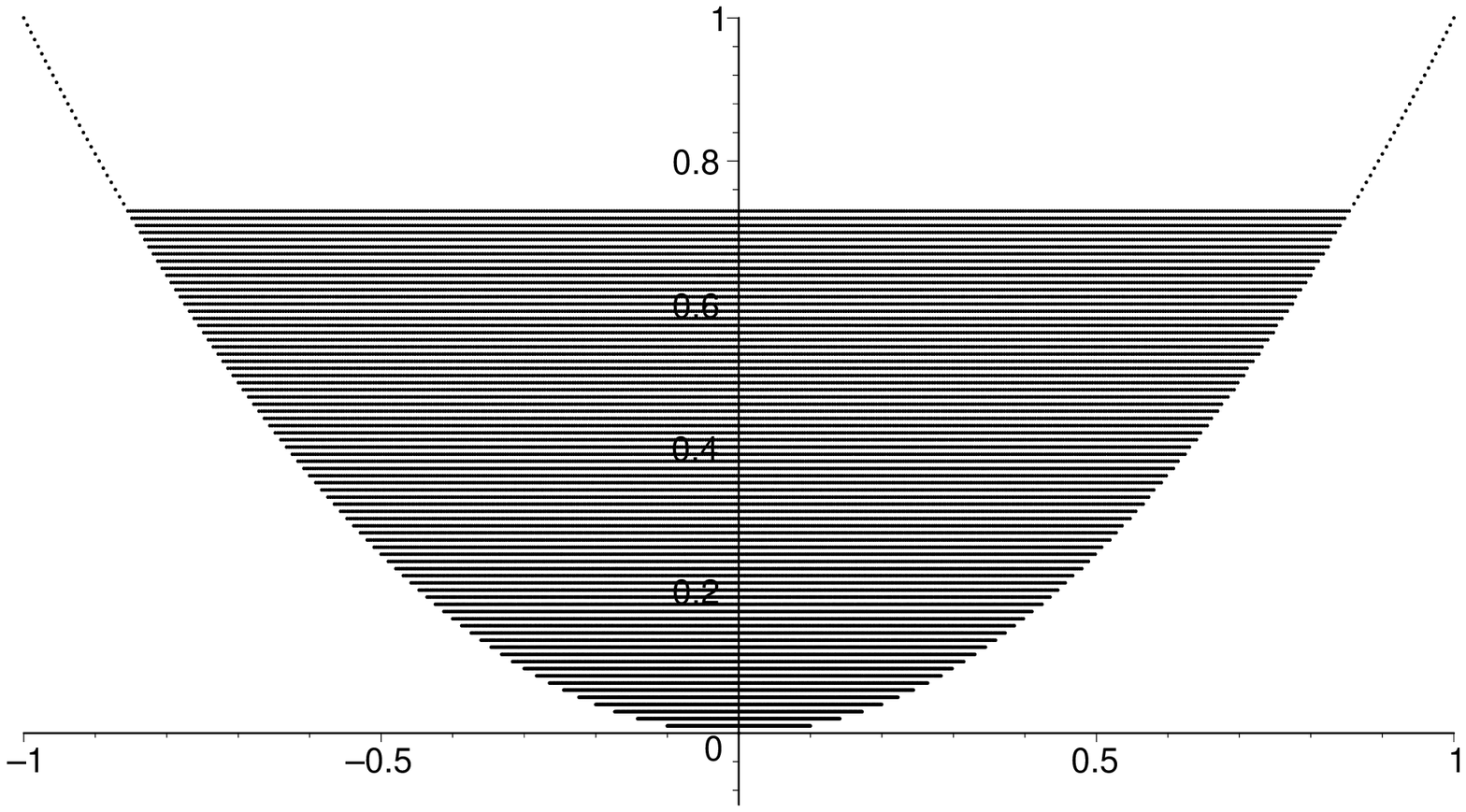}

\bigskip\bigskip
\noindent
{\sc Figure 3.} The dotted and shaded region is $\del\psi\big([0,1.5]\times\RR)$
for the Cauchy data $(\psi_0,\dot\psi_0)$ given by 
(\ref{PsizeroStwoEq}) and (\ref{PsiDotZeroEq}).

\vglue1cm

A typical set of the form $\del\psi\big([0,T]\times\RR^n\big)$ is illustrated in Figure 3.
It is contained in 
$\epi\,(-\dot u_0)\, \sm \epi\,(-(\dot u_0)^{\star\star})$, 
with $(\dot u_0)^{\star\star}(y)\equiv-1$.
The total curvature of $\h{graph}(\psi)$ is bounded from above by the area of the difference
between the epigraph of $-\dot u_0$ and that of $-(\dot u_0)^{\star\star}$ 
(over $[-1,1]$),
$$
\int_0^\infty\int_{-\infty}^\infty
\h{MA}\,\psi
=\int_{-1}^1(1-y^2)dy
=\frac43.
$$
Indeed, the total curvature is given by the measure in $S^n$ of
the unit normals to the surface, that is dominated from
above by the Lebesgue measure of its stereographic projection in $\RR^n$.
The latter is precisely the \MA mass (see \cite{RT}, Proposition 2.3, or \cite{P}).
Thus, while the surface has an infinitely long corner set, the surface quickly
becomes `almost linear' outside it, and in a manner guaranteeing
its total curvature will remain finite.

}
\end{exam}

\bigs
\section{Interior $C^1$ regularity of $\usdoublestar$}
\label{ConeRegularitySection}
\bigs

The purpose of this section is to prove Lemma \ref{LegendreLemma},
which shows that $\usd$ is $C^1$ on the interior of $P$.

Let $f\in C^0(P)$. The biconjugate $f^{\star\star}$ of $f$
is a convex function on $P$ and can be characterized as
the `convex envelope' of $f$ (\cite{HL2}, Theorem 1.3.5, p. 45)
\begin{equation}
\label
{BiconjugateFirstRepEq}
\begin{array}{lll}
f^{\star\star}(x)
& =\;\;\,
\sup
\big\{
g(x)\,:\, \h{\rm $g$ convex on $P$ and\ } g\le f
\big\}
\cr\cr
& =\;\;\,
\sup
\big\{
a(x)\,:\,  \h{\rm $a$ affine on $P$ and\ } a\le f
\big\}.
\end{array}
\end{equation}

\smallskip

Consider a smooth function defined
on a compact convex set and smooth in its interior.
It is not true in general that its biconjugate is differentiable
in the interior.
However, in our situation the biconjugate enjoys the maximal degree of regularity
possible in general, i.e., it is $C^{1,1}$ (no gain is achieved by considering
a real-analytic function).
General results in the literature (see \cite{BH,GR,KK} and \cite{HL2}, \S X.1.5)
are usually stated for functions defined on all of $\RR^n$ that obey certain
growth conditions at infinity (cf. \cite{GR}, (2.3),
\cite{HL2}, (1.5.2), p. 50, \cite{BH}, (18), \cite{KK}, (1))
or else make certain assumptions
regarding boundary behavior that need not hold in our situation (cf.
\cite{GR}, (4.1) and particularly Remark 4.3).
Therefore, and also since the constructions involved will be
useful later, we find it beneficial to extract from the references
above partially self-contained proofs of differentiability of
$\usd$. The $C^{1,1}$ estimate can be deduced from the $C^1$
estimate \cite{RZ3}.

\begin{lem}
\label{LegendreLemma}
Let $s>0$. 
\hfill\break
(i) The graph of the Legendre double dual of $u_s$ is the 
lower boundary of the closed convex hull of the
epigraph of $u_s$.
\hfill\break
(ii) $ u^{\star\star}_s\in C^{1}(P\setminus\partial P)\cap C^0(P)$.

\end{lem}

\begin{proof}
(i) From (\ref{BiconjugateFirstRepEq}) it follows that for any
finite convex function $f$ defined on $P$, one has
\begin{equation}
\label{CvxHullEpigraphEq}
\overline{\h{\rm co}\; \epi\, f} = \epi\, f^{\star\star},
\end{equation}
Since $u_s$ is smooth and continuous up to the boundary of $P$ the
convex hull of $\epi\, u_s$ is closed, and the result follows.

\smallskip
\noindent
(ii) From (\ref{BiconjugateFirstRepEq}) we have $\usd\le u_s$.
Since $u_s$ majorizes a linear function on $P$ it also
follows that $\usdoublestar$ is bounded below.
Hence $\usd\in C^0(P)$.

We divide the proof of the $C^1$ estimate
into two steps and closely follow \cite{BH,HL2}. The proof
will show that for any compact subset $\Omega\subset P\sm\del P$
there exists a compact subset $\Omega'\subset P\sm\del P$ with
$\Omega\subset \Omega'$ such that
$$
||\nabla\usd||_{C^0(\Omega)}
\le 
||\nabla u_s||_{C^0(\Omega')}
<C(s+1),
$$
with $C=C(\Omega,\Omega',u_0,\dot u_0)>0$.

\smallskip
\noindent
{\it First step.}
Denote by $\Delta_{n+1}\subset\RR^{n+1}$ the unit simplex
$$
\Delta_{n+1}:=\big\{\lambda=(\lambda_1,\ldots,\lambda_{n+1})\in\RR^{n+1}\,:\,
\lambda_i\ge0, \; \sum_{i=1}^{n+1}\lambda_i=1\big\}.
$$
Recall the following representation formula for the biconjugate function
(\cite{Ro}, Corollary 17.1.5),
\begin{equation}
\label
{BiconjugateRepEq}
u^{\star\star}_s(y)
=
\inf
\Big\{
\lambda\cdot(u_s(y_1),\ldots,u_s(y_{n+1}))
\,:\,
\lambda\in\Delta_{n+1}, \; y_i\in P,\; \sum_{i=1}^{n+1}\lambda_iy_i = y
\Big\}.
\end{equation}
We claim that there
exist---for each $y\in P\setminus\del P$---points
$y_1,\ldots,y_{n+1}\in P$ and a vector $\lambda\in\Delta_{n+1}$
such that
\begin{equation}
\label{CalledUponEq}
(y,\usdoublestar(y))=\sum_{i=1}^{n+1}\lambda_i(y_i,u_s(y_i)).
\end{equation}
To see that, observe that according to (i) the epigraph of $\usdoublestar$ is the closed
convex hull of the epigraph of $u_s$. Hence, by convexity, for each
$y\in P\setminus\del P$
there exists an $m\in\NN$ and a $\lambda\in\Delta_{m}$
and $\{(y_i,r_i)\}_{i=1}^m\subset \epi\; u_s$
such that
\begin{equation}
\label{UdoublestarConvexCombEq}
(y,\usdoublestar(y))=\sum_{i=1}^m\lambda_i(y_i,r_i).
\end{equation}
But then since $u_s(y_i)\le r_i$, and $\sum_{i=1}^m\lambda_i(y_i,u_s(y_i))$
also belongs to $\epi\, u_s$ it follows that whenever $\lambda_i>0$
there holds $r_i=u_s(y_i)$: otherwise $(y,\usdoublestar(y))$ would
lie ``directly above" another point in the convex hull of $\epi\; u_s$, contradicting
(i). Finally, using (i), since the convex hull of the epigraph of $u_s$ is closed, it follows
that $(y,\usdoublestar(y))$ lies on its boundary; then by
a consequence of Carath\'eodory's theorem it is
possible to take $m=n+1$ in (\ref{UdoublestarConvexCombEq}) (\cite{HL1}, Proposition 4.2.3, p. 126).

\smallskip
\noindent
{\it Second step.} When equation (\ref{CalledUponEq}) holds,
set
\begin{equation}
\label{ISetEq}
I:=\{\,i\,:\,\lambda_i>0\}.
\end{equation}
Following the terminology of \cite{BH}, \S3, we say
that the points $\{y_i\}_{i\in I}$ are {\it called upon by} $y$.
We omit from the notation the dependence of $y_i$ on $y$.
By (\ref{BiconjugateFirstRepEq}) $\usd$ is convex and $\usd\le u_s$.
Hence,
$$
\usdoublestar(y)
\le
\lambda\cdot(\usdoublestar(y_1),\ldots,\usdoublestar(y_{n+1}))
\le
\lambda\cdot(u_s(y_1),\ldots,u_s(y_{n+1})),
$$
which together with (\ref{CalledUponEq}) implies
\begin{equation}
\label{EqualityAtCalledUponPointsEq}
\usdoublestar(y_i)=u_s(y_i),\quad\forall i\in I.
\end{equation}
Note that whenever $g\le f$ and $g(x)=f(x)$  for $g$ convex and
$f$ differentiable,  then $g$ is differentiable at $x$ and $\nabla
g(x)=\nabla f(x)$ (indeed any element of $\del g(x)$ defines a
supporting hyperplane to $f$ at $x$, hence a tangent hyperplane).
\begin{claim}
\label{CalledUponInteriorClaim}
Let $y\in P\sm\del P$, and let $\{y_i\}$ and $I$ be defined by 
(\ref{CalledUponEq}) and (\ref{ISetEq}), respectively.
Then for each $i\in I$ one has $y_i\in P\setminus \del P$.
\end{claim}

\begin{proof}
The proof relies on the formula
\begin{equation}
\label{CalledUponSubdifferentialEq}
\del\usdoublestar(y)=\bigcap_{i\in I}\del u_s(y_i),
\end{equation}
(\cite{BH}, Theorem 3.6, \cite{HL2}, Theorem 1.5.6, p. 53)
whose derivation applies verbatim in our situation: indeed $y^*\in\del\usdoublestar(y)$
iff $u_s^{\star\star\star}(y^*)+\usdoublestar(y)=\langle y^*,y\rangle$, or
(using $u_s^{\star\star\star}=u_s^\star$ and (\ref{CalledUponEq}))
$\sum_{i\in I}\lambda_i(u_s^\star(y^*)+u_s(y_i)-\langle y^*,y_i\rangle)=0$,
and (\ref{CalledUponSubdifferentialEq}) follows by the characterization of the subdifferential.

Since $\usd\in C^0(P)$ and $y\in P\setminus \del P$, there exists a neighborhood
of $y$ in $P\setminus\del P$ on which $\usdoublestar$ is finite.
Hence, since $\usdoublestar$ is convex on $P$
it follows that the left hand side of (\ref{CalledUponSubdifferentialEq})
is nonempty (observe that this argument fails on the boundary: formally one thinks
of $u_s$ as defined on all of $\RR^n$ and identically equal to $+\infty$ outside $P$).
However, by Guillemin's formula (\ref{GuilleminFormulaEq})
the subdifferential of $u_G$ and hence
of $u_s$ at every point of $\del P$ is empty.
Therefore (\ref{CalledUponSubdifferentialEq}) implies that $y_i\in P\sm\del P$.
\end{proof}

Since $u_s\in C^1(P\setminus\del P)$
it follows from the Claim and (\ref{CalledUponSubdifferentialEq})
that $\usdoublestar$ is differentiable at $y_i,\; i\in I$,
and $\nabla\usdoublestar(y_i)=\nabla u_s(y_i)$
(observe that by (\ref{EqualityAtCalledUponPointsEq}) $y_i$ is the only
point called upon by $y_i$ and apply (\ref{CalledUponSubdifferentialEq})
to $y=y_i$).
Therefore, by using (\ref{CalledUponSubdifferentialEq}) once more,
it follows that $\del u_s(y_i)=\{\nabla u_s(y_i)\}$ is identical for all $i\in I$. We conclude
\begin{equation}
\label{CalledUponPointsGradientIdentityEq}
\nabla\usdoublestar(y)=\nabla\usdoublestar(y_i)=\nabla u_s(y_i),
\end{equation}
hence $\usdoublestar$ is differentiable at $y$, as desired.
Continuity of $\nabla\usdoublestar$ now follows from convexity (\cite{Ro}, Corollary 25.5.1,
\cite{KK}).
This concludes the proof of Lemma \ref{LegendreLemma}.
\end{proof}

\bigs
\section{The set $A_s$}
\label{AsSection}
\bigs

As indicated in \S\ref{IntroductionSection} the following 
set plays an important role in the analysis:

\begin{equation}
\label{AsDefEq}
A_s:=\{y\in P\,:\, u_s(y)\not=\usdoublestar(y)\}.
\end{equation}
In this section, we first prove a basic characterization of this
set which is later used repeatedly. In addition, we give
a geometric description of this set that is used to prove 
Lemma \ref{EssentialSmoothnessLemma} in the next
section. In \S\ref{MonotonicitySection} we prove further properties of $A_s$
which are needed for the proof of Theorem \ref{FirstMainThm}.

When $s\le T_\span^\cvx$ the set $A_s$ is empty.
When $s>T_\span^\cvx$ this set is an open non-empty
set relative to the topology of $P$ (it may intersect $\del P$---see 
Example \ref{AsBoundaryExample} below).
The set $A_s\cap (P\sm\del P)$ is a non-empty open set in the usual topology
of $\RR^n$. Both of these assertions follow from the continuity of $u_s$
and $\usdoublestar$.

\begin{lem}
\label
{AsAlternativeCharLemma}
One has
$$
A_s\;\cap\; P\sm\del P\,=\,\{\,y\in P\sm\del P\,:\, \del u_s(y)=\emptyset\}.
$$
\end{lem}

\begin{proof}
Let $y\in P\sm\del P$.
If $\del u_s(y)\not=\emptyset$ then $\usdoublestar(y)=u_s(y)$
(by definition $\usdoublestar$ majorizes the supporting affine function at $y$,
while $\usdoublestar\le u_s$ always),
hence
$$
A_s\;\cap\; P\sm\del P
\,\subset\,
\{\,y\in P\sm\del P\,:\, \del u_s(y)=\emptyset\}.
$$

Conversely, if $u_s(y)=\usdoublestar(y)$
then $y$ itself can be viewed as the only point called upon by $y$.
According to Lemma \ref{LegendreLemma}, convexity
of $\usdoublestar$, Equation (\ref{CalledUponSubdifferentialEq}),
and since $y\in P\sm\del P$, we have
\begin{equation}
\label
{GradientSubDiffBasicEq}
\del u_s(y)=\del\usdoublestar(y)
=\{\nabla u_s(y)\}
=\{\nabla\usdoublestar(y)\}
\ne\emptyset.
\end{equation}
Hence
$\{\,y\in P\sm\del P\,:\, \del u_s(y)=\emptyset\}
\,\subset\, A_s\;\cap\; P\sm\del P$.
\end{proof}

It is important to note that unlike the intuition from the case $n=1$ (see
\S\ref{OneDSection}),
the set $A_s$ may intersect the boundary of $P$. We illustrate with
a simple example.

\smallskip
\begin{exam}
\label{AsBoundaryExample}
{\rm 
We follow the notation of \S\ref{OneDSection}.
Consider $M=\PP^1\times\PP^1$ endowed with the Guillemin
\K structure $\o_{\vp_0}=\pi_1^*\oFS+\pi_2^*\oFS,$
where $\pi_j:M\ra \PP^1$ is the projection onto the $j$-th factor.
The associated polytope is $P=[-1,1]\times[-1,1]$, and
the symplectic potential dual to $\psi_0$ is
$$
u_0(y)=\half\sum_{j=1}^2(1+y_j)\log(1+y_j)+(1-y_j)\log(1-y_j).
$$
Let $f(y):\RR\ra\RR$ be a strictly concave function,
and let 
$$
\dot u_0(y)=f(y_1), \quad y\in P.
$$ 
Then for $s$ sufficiently large,
$$
A_s\cap \{\,(y_1,\pm 1)\,:\, y_1\in(-1,1)\,\}\ne\emptyset.
$$

}
\end{exam}
\medskip
\bigskip

Nevertheless, as illustrated in this example, what does generalize 
from the case $n=1$ is the following fact.

\begin{lem}
\label{AsPositiveDistLemma}
The set $A_s$ is at positive distance from the vertices of $P$.
\end{lem}

\begin{proof}
Recall that by the Delzant condition, 
if $p$ is a vertex of $P$
then $p$ is the intersection of exactly $n$ of the $d$ defining
half-spaces of $P$; in the notation of \cite{RZ2}, \S\S4.2,
\begin{equation}
\label{VertexPDefEq}
\{p\}=\bigcap_{k=1}^n\;\big\{\,y\in \RR^n\,:\,l_{j_k(p)}(y)=0\big\},
\end{equation}
with $j_k\in\{1,\ldots,d\}$. The functions
$l_{j_1(p)}(y),\ldots,l_{j_n(p)}(y)$ provide a coordinate system in which
the Hessian of $u_G$ is diagonal with eigenvalues
$\{\lambda_k=f_k l^{-1}_{j_k(p)}+g_k\}_{k=1}^n$, where
$f_k, g_k\in C^\infty(P)$ are bounded functions up to the boundary,
and $f_k>0$.
The same holds for $u_0\in\calL\calH(T)$ and hence also for $u_s$ since $u_s-u_0$ is smooth up
to the boundary.
Hence, for some $\eps>0$, that we assume is the largest such possible,
$u_s$ is strictly convex on $P\cap B(p,\eps)$, where 
$B(p,\eps):=\{v\in\RR^n\,:\, |v-p|<\eps\}$.

Let now $\delta\in(0,\epsilon)$ and
let $y\in   B(p,\delta)\cap\, P\sm\del P$. 
By convexity the graph of the tangent hyperplane
$H(w):=u_s(y)+\langle w-y,\nabla u_s(y)\rangle$
to the graph of $u_s$ at $y$ supports the graph 
of $u_s$ above $P\cap B(p,\eps)$. Suppose that this 
is not true globally, namely, 
$\min_P u_s-H <0$. From the explicit expression for $u_G$
we obtain that in the aforementioned coordinates the gradient
of $u_s$ can be written in the form 
$(\log l_{j_1(p)}+h_1,\ldots,\log l_{j_n(p)}+h_n)$,
where $h_j\in C^\infty(P),\,j=1,\ldots,n$. Since $u_s$ is bounded on $P$, 
it follows that by making $\delta>0$ smaller (and hence
making $y$ closer to $p$), we may assume that
for some $w\in P\sm(\{y\}\cup\del P)$, the graph of $H$ 
is tangent to that of $u_s$ at $w$, namely 
$\nabla u_s(y)=\nabla u_s(w)$. However, 
we already know that $w\not\in B(p,\eps)$, and
by evaluating the expression for $\nabla u_s$ 
with respect to the coordinates $l_{j_1(p)},\ldots,
l_{j_n(p)}$ at $w$ and comparing with (\ref{VertexPDefEq})
we obtain a contradiction
once $\delta$ is chosen sufficiently small with
respect to $\epsilon$ (in a manner depending on
$\max_j\max_P h_j$). It follows that for
such a choice of $\delta>0$ (that depends only on $p$
and $s$) we have 
$(B(p,\delta)\cap P)\cap A_s=\emptyset$.
\end{proof}

\bigs

\section{Essential smoothness of $\usdoublestar$}
\label{EssentialSmoothnessSection}

\bigs

Having established the interior $C^1$ regularity of
$\usdoublestar$ we now turn to a result concerning the boundary
behavior of its gradient.

\begin{lem}
\label{EssentialSmoothnessLemma}
For each $s>0$ the function
$\usdoublestar$ is essentially smooth.

\end{lem}

\begin{proof}
First observe that by Lemma \ref{LegendreLemma} (i) the function
$\usdoublestar$ is proper, as required by Definition \ref{ConvexTermsDef}.
Next, observe that $\usdoublestar$ is differentiable
on $P\sm\del P=\h{\rm int}\,\dom(\partial \usdoublestar))$ by Lemma
\ref{LegendreLemma} (ii).
Also, from Lemmas \ref{LegendreLemma} and \ref{AsPositiveDistLemma}
and their proofs,
if $\{w_i\}\subset P\sm\del P$ is a sequence converging
to one of the vertices of $P$ then
$\lim_{i\ra\infty}|\nabla \usdoublestar(w_i)|=+\infty$
(since $\nabla\usdoublestar(w_i)=\nabla u_s(w_i)$ for
large $i$).
Consider now a sequence $\{w_i\}\subset P\sm\del P$ converging
to a point $p$ in $\del P$ contained in the interior of
a face $F$ of dimension $1$, that we assume,
without loss generality, is cut out
by the equations $l_j=0,\;$ for $j=1,\ldots,n-1$, with
$l_n\in[0,C]$ a coordinate on this face.
Then $\lim_{i\ra\infty} l_j(w_i)=0$ for $j=1,\ldots,n-1$.
Let $\{p_i^0\},\{p_i^1\}\subset
P\sm\del P$ be sequences defined by the equations
$l_k(p_i^j)=l_k(w_i), k=1,\ldots,n-1$ and
$l_n(p_i^j)=Cj+(1/2-j)\epsilon/i$, with $j=0,1$, where $\epsilon>0$
is some sufficiently small constant depending on $p$ and $P$.
Then $\lim_{i\ra\infty} p_i^j,$ with $j=0,1,$ are two distinct vertices of $P$. Using the functions
$l_1,\ldots,l_n$ as coordinates in a neighborhood of the face $F\subset P$,
it follows from Guillemin's formula (\ref{GuilleminFormulaEq})
that for each $k\in\{1,\ldots,n-1\}$ one has
$$
L_k:=\lim_{i\ra\infty}\frac{\del \usdoublestar}{\del l_k}(p_i^0)=
\lim_{i\ra\infty}\frac{\del \usdoublestar}{\del l_k}(p_i^1)\in\{\pm\infty\},
$$
since by Lemma \ref{AsPositiveDistLemma} one may replace $\usdoublestar$
by $u_s$ in this equation.
Considering then the function $\usdoublestar$ restricted
to the line connecting $p_i^0$ and $p_i^1$ (necessarily contained in $P$ by
convexity), and taking the limit
it follows that
$$
\lim_{i\ra\infty}\frac{\del \usdoublestar}{\del l_k}(w_i)=L_k,
$$
and so $\lim_{i\ra\infty}|\nabla \usdoublestar(w_i)|=+\infty$.
The general case where $p$ is contained in a boundary face of dimension
$m\le n-1$ now follows by induction on $m$, using arguments as above.
\end{proof}

\bigs

\section{Gradient and subdifferential mappings of $\usdoublestar$}
\label{GradientSubDiffUsSection}

\bigs
The surjectivity of $\nabla\usd$ can be proved as follows.
By definition,
$$
\psi_s(x)=\sup_{y\in P}\big(\langle x,y\rangle -u_s(y) \big),
$$
and since $P$ is compact and $u_s$ bounded it follows that the
supremum is achieved at some $y\in P$. Duality
then implies that $x\in\del u_s(y)=\del\usd(y)$
(\cite{HL2}, Theorem 1.4.1, p. 47), and essential
smoothness (Lemma \ref{EssentialSmoothnessLemma}) 
implies that $y\in P\sm\del P$. Hence by
Lemma \ref{LegendreLemma} we have
$\del u_s(P\sm\del P)=\nabla \usd(P\sm\del P)=\RR^n$,
as desired.

Our goal in this section is to prove the surjectivity of $\nabla\usd$
by using an alternative homotopy argument. 
An advantage of this approach is that 
in the course of the proof we also obtain that $\usd$ and 
$\nabla\usd$ are continuous in $s$.  This fact is useful
when studying the structure of the singular locus of $\psi$
(see the discussion in \S\ref{IntroductionSection} just
before \S\S\ref{QuantumConsequencesSubsection}).
We believe that albeit being more involved, this approach
has its own interest, and might also find applications in
situation where the families of functions studied are less
explicit. 

First, we state an elementary lemma that describes the gradient
image of $u_s$.

\begin{lem}
\label
{SurjectiveGradientImageLemm}
For each $s\in\RR_+$, we have
$$
\nabla u_s(P\setminus\partial P)
=
\RR^n.
$$
\end{lem}

\begin{proof}
Recall that by (\ref{LHTDefSecondEq})
we have
\begin{equation}
\label{UzeroUGDifferenceEq}
u_0-u_G\in C^\infty(P).
\end{equation}
From the explicit formula
(\ref{GuilleminFormulaEq}) for $u_G$ and the Delzant conditions
one may verify directly that $ \nabla u_G(P\setminus\partial P)
=\RR^n. $ 
Note that (\ref{GuilleminFormulaEq}) and (\ref{UzeroUGDifferenceEq}) also imply
that $\nabla u_s$ is, as a smooth map of $P\sm\del P$ into 
$\RR^n$, properly homotopic to 
$\nabla u_0$, that is itself properly homotopic to $\nabla u_G$.
It
follows that the topological degree of $\nabla u_s$ equals that of
$\nabla u_G$. Since $\nabla u_G$ is a bijection we have $\h{\rm
deg}\,(\nabla u_s)=1$. It follows that $\nabla u_s:P\sm\del P\ra
\RR^n$ is surjective (see, e.g., \cite{OR}, Chapter 3, or
\cite{G2}, Theorems 3.6.6, 3.6.8).
\end{proof}

\begin{lem}
\label
{SurjectivityRnLemma}
Let $s>0$.
\hfill\break
(i) We have,
$$
\nabla u^{\star\star}_s({P\setminus\partial P})=\del u_s(P\sm\partial P)=\RR^n.
$$
\hfill\break
(ii) Moreover,
\begin{equation}
\label{SubDiffImageUsAsEq}
\del u_s(P\sm(\partial P\cup A_s))
=
\RR^n,
\end{equation}
and
$$
\nabla\usdoublestar=\nabla u_s=\del u_s,\quad \h{\ on\ } \quad P\sm (\del P\cup A_s).
$$

\end{lem}

\begin{proof}
\smallskip
\noindent
(i)
Equation (\ref{CalledUponSubdifferentialEq})
implies that $\del\usdoublestar(P\sm\del P)\subset\del u_s(P\sm\del P)$. 
On the other hand, if $\del u_s(y)\not=\emptyset$ then
by Lemma \ref{AsAlternativeCharLemma} and its proof
we have $\usdoublestar(y)=u_s(y)$ and
$\del\usdoublestar(y)=\del u_s(y)$,
hence $\del u_s(P\sm\del P)\subset \del \usdoublestar(P\sm\del P)$.
Thus,
\begin{equation}
\label{SubdiffImageUsdoublestarEq}
\del u_s(P\sm\del P)=\del \usdoublestar(P\sm\del P).
\end{equation}
Since (by Lemma \ref{LegendreLemma})
$\nabla\usdoublestar(P\sm\del P)=\del\usdoublestar(P\sm\del P)$,
it suffices to show
\begin{equation}
\label{SubdiffImageUsEq}
\nabla \usdoublestar(P\sm\del P)=\RR^n.
\end{equation}

To that end, we will again rely on degree theory, however
now for continuous maps of $S^n$.
We extend $\nabla\usdoublestar$ to a map
$G(s,y):S^n\ra S^n$, where $S^n=\RR^n\cup\{\infty\}$, defined
by
$$
G(s,y)=
\begin{cases}
\nabla\usdoublestar(y),&  y\in P,
\cr\cr
\infty,&  y\in S^n\sm P.
\end{cases}
$$
Observe that
by Lemma \ref{EssentialSmoothnessLemma} for each fixed $s$
the map $G(s,\,\cdot\,):S^n\ra S^n$ is continuous, and
\begin{equation}
\label
{GsyEq}
G(s,y)=
\begin{cases}
\nabla\usdoublestar(y),& y\in P\sm\del P,
\cr\cr
\infty,&  y\in \overline{S^n\sm P}.
\end{cases}
\end{equation}
Next, we prove that this map is also continuous in $s$.

\begin{claim}
\label
{HomotopyClaim}
As continuous maps, $\nabla\usdoublestar$ is homotopic to $\nabla u_0$.
\end{claim}

\begin{proof}
By definition, we need to show that the map
$G:\RR_+\times S^n\ra S^n$ defined above is continuous.
As pointed out, it remains to show that for each fixed
$y\in S^n$, the map $G(\,\cdot\,,y):\RR_+\ra S^n$ is continuous.
Observe that by (\ref{GsyEq}) it only remains to treat
the case $y\in P\sm \del P$.

Fix $y\in P\sm\del P$ as well as $s>0$.
Let $\{s_j\}_{j\ge1}\subset\RR_+$ satisfy $\lim_{j\ra\infty}s_j=s$.
By convexity and Lemma \ref{LegendreLemma} (ii) we have for each $j\ge1$,
\begin{equation}
\label
{UsdSubDiffSeqEq}
u_{s_j}^{\star\star}(y')
\ge
u_{s_j}^{\star\star}(y)
+
\langle\, y'-y,\nabla u_{s_j}^{\star\star}(y)\,\rangle,
\quad \forall y'\in P.
\end{equation}
Let $y^*$ be any limit point of
$\{ \nabla u_{s_j}^{\star\star}(y) \}_{j\ge1}$.
If for each $v\in P$ the map $s\mapsto \usd(v)$ were continuous
then taking the limit as $j$ tends to infinity in
(\ref{UsdSubDiffSeqEq}) would imply
$$
\usd(y')\ge \usd(y)
+
\langle\, y'-y,y^*\,\rangle,
\quad \forall y'\in P.
$$
implying that $y^*\in\del\usd(y)$.
Lemma \ref{LegendreLemma} would then imply
that $y^*=\nabla\usd(y)$, proving the Claim.
Now, to prove the continuity of $\usd(v)$ in $s$ it suffices
to use the representation formula (\ref{BiconjugateRepEq}).
On the one hand,

\begin{equation*}
\begin{array}{lll}
u^{\star\star}_{s'}(y)
&\! = &\dis
\!\!\inf
\Big\{
\sum_{i=1}^{n+1}\lambda_iu_{s'}(y_i)
\,:\,
\lambda\in\Delta_{n+1}, \; y_i\in P,\; \sum_{i=1}^{n+1}\lambda_iy_i = y
\Big\}
\cr
& \! = &\dis
\!\!\inf 
\Big\{
\sum_{i=1}^{n+1}\lambda_iu_{s}(y_i)+(s'-s)\!\sum_{i=1}^{n+1}\lambda_i\dot u_0(y_i)
\!:
\lambda\in\mskip-1mu\Delta_{n+1}, y_i\in\mskip-1mu P, \sum_{i=1}^{n+1}\lambda_iy_i = y\!
\Big\}
\cr
& \!\ge &\dis
\!\!\inf
\Big\{
\sum_{i=1}^{n+1}\lambda_iu_{s}(y_i)
\,:\,
\lambda\in\Delta_{n+1}, \; y_i\in P,\; \sum_{i=1}^{n+1}\lambda_iy_i = y
\Big\}
\cr
&& \dis + 
\inf
\Big\{
(s'-s)\sum_{i=1}^{n+1}\lambda_i\dot u_0(y_i)
\,:\,
\lambda\in\Delta_{n+1}, \; y_i\in P,\; \sum_{i=1}^{n+1}\lambda_iy_i = y
\Big\}
\cr
& \! \ge & \dis
\!\!\usd(y)+(s'-s)\min_P\dot u_0.
\end{array}
\end{equation*}
And on the other hand,

\begin{equation*}
\begin{array}{lll}
u^{\star\star}_{s'}(y)
& \le &\dis
\inf 
\Big\{
\sum_{i=1}^{n+1}\lambda_iu_{s}(y_i)+\displaystyle (s'-s)\max_P\dot u_0
\,:\,
\lambda\in\Delta_{n+1}, \; y_i\in P,\; \sum_{i=1}^{n+1}\lambda_iy_i = y
\Big\}
\cr
& = &\dis
\usd(y)+(s'-s)\max_P\dot u_0.
\end{array}
\end{equation*}
Hence,
$$
|u^{\star\star}_{s'}(v)-u^{\star\star}_s(v)|
<
(s'-s)\cdot||\dot u_0||_{C^0(P)}<C(s'-s),
$$
as desired.
\end{proof}

Thus, applying Lemma \ref{SurjectiveGradientImageLemm} for $s<T_\span^\cvx$,
it follows that
$\deg\,\nabla\usdoublestar=1$, and hence that $\nabla\usdoublestar$ is
surjective (see, e.g., \cite{H}, \S2.2). This implies (\ref{SubdiffImageUsEq}).

\smallskip
\noindent (ii) The second part of the statement is already
contained in Lemma \ref{AsAlternativeCharLemma} and its proof (see
(\ref{GradientSubDiffBasicEq})). The first part follows from (i)
and Lemma \ref{AsAlternativeCharLemma}.
\end{proof}

\bigs
\section{Monotonicity and continuity of $A_s$}
\label{MonotonicitySection} 
\bigs

The following two lemmas establish the continuity of
the set-valued map $s\mapsto A_s$, its strict 
monotonicity, and identify its asymptotic limit
$A_\infty$. Later, in Lemma \ref{MASingMassAinfinityLemma}, 
the first two facts are used to establish
a strict lower bound on the 
\MA  mass of the Legendre transform potential,
while the third fact is used to show an a priori
upper bound.

\begin{lem}
\label{AinfinityLemma}
Let $s_2>s_1>0$. Then
\begin{equation}
\label{AsTimeInclusionEq}
\overline{A_{s_1}}\sm\del P\subset \h{\rm int}\, A_{s_2}.
\end{equation}
In addition,
\begin{equation}
\label{AinfinityEq}
A_\infty:=
\bigcup_{s>0}A_s
=\{y\in P\,:\, \dot u_0(y)\ne (\dot u_0)^{\star\star}(y)\},
\end{equation}

\end{lem}

\begin{proof}
We claim that whenever $s_1<s_2$ one has
\begin{equation}
\label{AsTimeInclusionSecondEq}
A_{s_1}\subset A_{s_2}.
\end{equation}
Let $y\in A_{s_1}\sm\del P$. By Lemma \ref{AsAlternativeCharLemma}
then $\del u_{s_1}(y)=\emptyset$.
Hence,
$\nabla u_{s_1}(y)$ is not a subgradient for $u_{s_1}$ at $y$,
so there exists $y'\in P\sm\{y\}$ such that
\begin{equation}
\label{UsAboveGraphEq}
u_{s_1}(y')<\langle\, y'-y,\nabla u_{s_1}(y)\,\rangle+u_{s_1}(y).
\end{equation}

\smallskip
\noindent
On the other hand, since $u_0$ is strictly convex, 
\begin{equation}
\label{UzeroStrictConvexityIneqEq}
u_0(y')>\langle\, y'-y,\nabla u_0(y)\,\rangle+u_0(y),
\end{equation}
and
it follows that
\begin{equation}
\label{DotUzeroAboveGraphEq}
\dot u_0(y')<\langle\, y'-y,\nabla \dot u_0(y)\,\rangle+\dot u_0(y).
\end{equation}
Multiplying (\ref{DotUzeroAboveGraphEq}) by $s_2-s_1>0$ and adding to
(\ref{UsAboveGraphEq}) then implies that $\del u_{s_2}(y)=~\emptyset$,
proving (\ref{AsTimeInclusionSecondEq}).

To prove (\ref{AsTimeInclusionEq}) it remains to
show that $\del A_{s_1}\subset\h{\rm int}\,A_{s_2}$.
Assume that $y\in\del A_{s_1}$. Then there exists some $y'\in P\sm\{y\}$
such that
\begin{equation}
\label{USoneEqualityEq}
u_{s_1}(y')=\langle\, y'-y,\nabla u_{s_1}(y)\,\rangle+u_{s_1}(y).
\end{equation}
It then follows from (\ref{UzeroStrictConvexityIneqEq}) that
(\ref{DotUzeroAboveGraphEq}) holds.
As before, this implies that
$$
u_{s_2}(y')<\langle\, y'-y,\nabla u_{s_2}(y)\,\rangle+u_{s_2}(y),
$$
hence $y\in\h{\rm int}\,A_{s_2}$, as desired.

Next, note that (\ref{DotUzeroAboveGraphEq})
implies that $\del \dot u_0(y)=\emptyset$, hence
$$
A_s\subset\{y\in P\,:\, \dot u_0(y)\ne (\dot u_0)^{\star\star}(y)\},
$$
for every $s>0$.

Conversely, given $y\in P$ such that
$\del \dot u_0(y)=\emptyset$, let $y'\in P\sm\del P$ be such that
(\ref{DotUzeroAboveGraphEq}) holds. It follows that (\ref{UsAboveGraphEq})
must hold for $s_1>0$ large enough, i.e., $\del u_{s_1}(y)=\emptyset$
and $y\in A_{s_1}$. By (\ref{AsTimeInclusionEq})
then $y\in A_s$ for every $s>s_1$ and (\ref{AinfinityEq}) follows.
\end{proof}

\begin{lem}
\label{AsContinuousMappingLemma}
The map $s\mapsto A_s$ is continuous as a set-valued map.
\end{lem}

\smallskip

\begin{proof}
We will show that the map is both lower and upper semi-continuous.
Namely, for given $s$ and $\epsilon>0$ there exists $\delta>0$ such that for
all $s'\in (s-\delta,s+\delta)$ holds
\begin{equation}
\label{LSCSetValuedEq}
A_s\subset A_{s'}+B(0,\epsilon),
\end{equation}
and
\begin{equation}
\label{USCSetValuedEq}
A_{s'}\subset A_s+B(0,\epsilon),
\end{equation}
where $B(0,\epsilon):=\{y\in\RR^n\,:\, |y|<\epsilon\}$,
and the sum is in the sense of Minkowski.

First, we prove the lower semi-continuity. 
Let $y\in A_s$.
We start with the special case where
$\delta_s:=d(A_{s},\del P)>0$, where $d$
denotes the Euclidean distance function.
Hence,
by the preceeding Lemma we may assume that $s'<s$.
Consider the function $F:\RR_+\times (P\sm\del P)\times P\ra\RR$ defined by
$$
F(\sigma,v,w):=u_\sigma(w)-u_\sigma(v)-\langle w-v, \nabla u_\sigma(v)\rangle.
$$
Note that $F$ is smooth on its domain.
Let $G:\RR_+\times P\sm\del P\ra\RR$ be defined by
$$
G(\sigma,v):=\min_{w\in P}F(\sigma,v,w).
$$
Then $G$ is continuous. Since $F$ is uniformly Lipschitz
on $[s-1,s]\times P\sm B(\del P,\delta_s)\times P$, it follows
that $G$ is uniformly continuous on $[s-1,s]\times P\sm B(\del P,\delta_s)$
(we assume without loss of generality that $s>1$),
where 
$B(\del P,\delta_s)$ denotes a $\delta_s$-neighborhood of $\del P$
in $P$.
Note that in general $A_\sigma\sm\del P=\{y\in P\sm\del P\,:\, G(\sigma,y)<0\}$.
By the preceeding Lemma we have $d(A_{s'},\del P)\ge\delta_s$.
Hence, under our assumption that $\delta_s>0$, 
it follows that there exists some $\delta>0$
such that (\ref{LSCSetValuedEq}) holds whenever $|s'-s|<\delta$.

We now turn to the general case, and 
assume $A_s\cap \del P\ne\emptyset$.
By the previous case, we already know that 
we may choose $\delta=\delta(\delta_s)>0$ 
such that for every $y\in A_s$ with $d(y,\del P)\ge\delta_s$ 
there exists some $y'\in A_{s'}$ satisfying $|y-y'|<\eps$.
However, we need to show that $\delta$ does not tend to zero
as $\delta_s$ does.

To that end, let us assume that $y$ is close
to the boundary of $P$, but not in $\del P$.
We will return to the case $y\in \del P$ later.
We assume also, without loss of generality,
that $l_1(y)=\min_{i\in\{1,\ldots,d\}}l_j(y)$,
with $l_1(y)<l_j(y)$ for all $j\in\{2,\ldots,d\}$.
We complement $l_1$ with $n-1$ other functions,
which for simplicity of notation we assume
are $l_2,\ldots,l_n$, in such a manner that 
$l_1,\ldots,l_n$ form a coordinate system in $\RR^n$.

Now, we consider the function $F$ with its last 
argument restricted to the line 
$L:=\{v\in P\,:\, l_j(v)=l_j(y),\,j=2,\ldots,n\}$ in $P$ 
that passes through $y$ and is perpendicular to the face 
$$
\calF:=\{v\,:\,l_1(v)=\langle v,v_1\rangle-\lambda_1=0\},
$$
i.e., 
$$
H(\sigma,y,t):=u_\sigma(y+tv_1)-u_\sigma(y)-\langle tv_1, \nabla u_\sigma(y)\rangle,
\quad t\in[-C_1,C_2],
$$
with $C_1=C_1(y)>0$ proportional to $d(y,\del P)$ (up to some uniform constant),
and such that $l_1(y-C_1v_1)=0$,
and with $C_2=C_2(y)>0$ uniformly bounded from above and away from zero,
(under the assumption that $l_1(y)<\delta_s$), and such that
$y+C_2v_1\in \del P$.

The last term in $H$, computed with respect to the 
coordinates $l_1,\ldots,l_n$, equals
$$
\langle w-v, \nabla u_\sigma(v)\rangle
=
a_\sigma(t)\log l_1(y)+b_\sigma(t),\qquad (\hbox{where \ } w=v+tv_1),
$$
for some uniformly bounded functions $a_\sigma(t),b_\sigma(t)$ of $t\in[-C_1,C_2]$,
and with $a_\sigma(t)<0$ for $t>0$, and $a_\sigma(t)\ge0$ for $t\le0$.
Hence, 
$$
H(\sigma,y,t)<C'-a_\sigma(t)\log l_1(y),
$$
for some $C'=C'(\sigma)>0$. 
It then follows that if $y$
is taken close enough to $\calF$, i.e., $-\log l_1(y)$
is large enough, then for some some $t\in[0,C_2)$
we will have $H(\sigma,y,t)<0$.

Note that $C'=C'(\sigma)$ is uniformly bounded
for $\sigma$ in some small neighborhood in $\RR_+$ 
(independently of $y$). Hence, the above arguments
imply that given $\epsilon>0$, we may find $\delta>0$
such that whenever $|s'-s|<\delta$ we may also
find $y'\in P\sm\del P$ with $|y'-y|<\eps$ and such that 
$H(s',y',t')<0$ with $t'\in(0,C_2(y'))$.
In particular, we found a $w':=y'+t'v_1$ satisfying
$F(s',y',w')<0$, hence $y'\in A_{s'}$, as desired.

Finally, if $y\in \del P$, since $A_s$ is open
in the relative topology of $P$ we may choose
$\tilde y\in A_s\cap P\sm\del P$ with 
$|\tilde y-y|<\eps/2$. We may then carry out the arguments
above for $\tilde y$ to find $\delta>0$ such that whenever
$|s'-s|<\delta$ there exists $y'\in A_{s'}$ with
$|y'-\tilde y|<\eps/2$. Hence, once again, 
$y\in A_{s'}+B(0,\eps)$, and this concludes the proof
of the lower semi-continuity.

The proof of the upper semi-continuity (\ref{USCSetValuedEq}) 
involves similar arguments, by switching the roles of $s$ and $s'$.
\end{proof}

\bigs
\section{On the invertibility of $\nabla \usdoublestar$ and strict
convexity of $\usd$}
\label
{InvertibilitySection}
\bigs

In this section we describe the set on which $\usd$ is strictly convex.

\begin{lem}
\label{AsInvertibilityCharLemma}
$\nabla \usdoublestar$ is invertible on $\h{\rm int}\,\big(P\setminus(\partial P\cup A_s)\big)$.
Moreover, if $y\in \h{\rm int}\,\big(P\setminus(\partial P\cup A_s)\big)$
then
\begin{equation}
\label{StrictConvexityUsdoublestarSupportFirstEq}
u_s(y')
>\langle\, y'-y,\nabla u_s(y)\,\rangle+u_s(y),\quad \forall y'\in P\sm\{y\},
\end{equation}
and
\begin{equation}
\label{StrictConvexityUsdoublestarSupportSecondEq}
\usdoublestar(y')>\langle\, y'-y,\nabla\usdoublestar(y)\,\rangle+\usdoublestar(y),\quad \forall
y'\in P\sm\{y\}.
\end{equation}
\end{lem}

\begin{proof}
Note that (\ref{StrictConvexityUsdoublestarSupportSecondEq})
implies (\ref{StrictConvexityUsdoublestarSupportFirstEq}),
since $\usdoublestar(y')\le u_s(y')$, $\usdoublestar(y)= u_s(y)$,
and $\nabla\usdoublestar(y)=\nabla u_s(y)$. In addition, the
invertibility statement also follows from
(\ref{StrictConvexityUsdoublestarSupportSecondEq}).

Let $y\in\h{\rm int}\,\big(P\setminus(\partial P\cup A_s)\big)$.
By convexity of $\usdoublestar$,
\begin{equation}
\label{StrictConvexityUsdoublestarSupportThirdEq}
\usdoublestar(y')\ge\langle\, y'-y,\nabla\usdoublestar(y)\,\rangle+\usdoublestar(y),
\quad\forall y'\in P.
\end{equation}
Suppose now that equality holds in
(\ref{StrictConvexityUsdoublestarSupportThirdEq})
for some $y'\in P$ (by Lemma \ref{EssentialSmoothnessLemma}
necessarily $y'\in P\sm\del P$). By convexity
and differentiability of $\usdoublestar$ then
\begin{equation}
\label{StrictConvexityUsdoublestarSupportFourthEq}
\nabla\usdoublestar(y')=\nabla\usdoublestar(y)=\nabla u_s(y).
\end{equation}
Since $\del u_s(y)\ne\emptyset$ then $u_s$ is convex at $y$,
i.e., $\nabla^2u_s(y)\ge0$. We claim that in fact
\begin{equation}
\label{StrictConvexityUsAtIntPsmAsEq}
\nabla^2 u_s(y)>0.
\end{equation}
Before proving (\ref{StrictConvexityUsAtIntPsmAsEq})
let us prove that it implies that $y'=y$, and hence
proves (\ref{StrictConvexityUsdoublestarSupportSecondEq}).
Indeed, consider the function
$$
F(t):=\usdoublestar(ty+(1-t)y').
$$
This function is convex for $t\in[0,1]$, and
strictly convex for $t\in(1-\epsilon,1]$ for some $\epsilon\in(0,1)$
(here we use (\ref{StrictConvexityUsAtIntPsmAsEq}) and that
$F(t)=u_s(ty+(1-t)y')$ for $t$ near $1$ since then
$ty+(1-t)y'\in\h{\rm int}\,\big(P\setminus(\partial P\cup A_s)\big)$).
Thus,
$$
\langle\nabla\usdoublestar(y'),y-y'\rangle
=
F'(0)<F'(1)
=
\langle\nabla\usdoublestar(y),y-y'\rangle,
$$
By (\ref{StrictConvexityUsdoublestarSupportFourthEq})
we must have then $y'=y$.

We return to proving (\ref{StrictConvexityUsAtIntPsmAsEq}).
Assume on the contrary that $\langle \nabla^2 u_s(y)\xi,\xi\rangle=0$
for some $\xi\in\RR^n\sm\{0\}$.
First, observe that since $\nabla^2 u_0>0$,
it follows that for every $\epsilon>0$ one has
$\langle\nabla^2 u_{s+\epsilon}(y)\xi,\xi\rangle<0$. Hence,
$\del u_{s+\epsilon}(y)=\emptyset$, i.e., $y\in A_{s+\epsilon}\sm\del P$.
By Lemma \ref{AinfinityLemma} it follows that $y\in\overline{A_s}\sm\del P$,
contradicting our assumption.
\end{proof}

A certain converse of the preceeding Lemma is given by the following.

\begin{lem}
\label
{WhereUsdStrictlyCvxLemma}
(i) Let $y\in P\sm\del P$. Assume that $\usdoublestar$ is strictly convex
in some neighborhood of $y$.
Then $y\in \h{\rm int}\,\big(P\sm(\del P\cup A_s)\big)$.
\hfill\break
\noindent
(ii) Let $y\in A_s\sm\del P$. Then there exists a line 
in $\overline{A_s}\sm\del P$ passing through $y$ 
and intersecting $\del A_s\sm\del P$, along 
which $\nabla\usdoublestar$ is constant.

\end{lem}

\begin{proof}
We only prove (ii) since it implies (i). However, note that
(i) is also a consequence of Lemma \ref{SurjectivityRnLemma} (ii)
that implies that $\del\usd(A_s\sm\del P)=\emptyset$.

Let $y\in A_s\sm\del P$.  Consider the tangent
hyperplane at $y$, given by the equation $l(y')=\usdoublestar(y)+\langle
\nabla\usdoublestar(y),y'-y\rangle$. By convexity $l\le
\usdoublestar$ on $P$. If one has
$l<\usdoublestar=u_s$ on $P\sm  A_s$,
then by compactness for some $\epsilon>0$
one has also has $l+\epsilon<\usdoublestar$ there. However,
by (\ref{BiconjugateFirstRepEq})
a fortiori $\usd$ equals the supremum of all affine
functions majorized by $u_s$ over $P\sm A_s$.
But then one would obtain a contradiction
to $l(y)=\usdoublestar(y)$. It follows that for some
$y'\in P\sm A_s$ one has $l(y')=\usdoublestar(y')$. Note that
by the essential smoothness of $\usdoublestar$ proved in
Lemma \ref{EssentialSmoothnessLemma} we must have
$y'\in P\sm \del P$.
Since $l$ is affine, convexity then implies that for each
point on the line segment connecting $y$ to $y'$ one has $l=\usdoublestar$.
It follows that $l$ is the tangent hyperplane to $\usdoublestar$
at each one of those points. Since $y'\not\in\h{\rm int}\, A_s$
it follows that the line connecting $y$ and $y'$ intersects $\del A_s\sm\del P$,
proving our claim.

As an alternative proof, one may also show that $\usd$ is affine on the polyhedron
$\h{\rm co}\,\{y_i\,:\,~i\in~I\}$ containing $y$ and generated by the points called upon by
$y$ (\cite{HL2}, Theorem 1.5.5, p. 52).
\end{proof}

Note that the preceeding Lemma does not quite give
a foliation of $A_s\sm\del P$ by lines along which
$\nabla \usdoublestar$ is constant. The
rank of $\ker\nabla^2\usdoublestar$ may jump in $A_s\sm\del P$
and so there may be more than one line with that property
passing through a given point. Instead, $A_s\sm\del P$ is partitioned
into maximal sets along which $\nabla \usdoublestar$ is constant
(see (\ref{QMaxSYDefEq}), (\ref{AsPartitionEq}) and \S\ref{MAMeasureLegendrePotentialSection}).

Yet, as a corollary of Lemma \ref{WhereUsdStrictlyCvxLemma} (ii)
(or of Lemma \ref{SurjectivityRnLemma} (ii)) 
we have
$$
\del \usdoublestar(A_s\sm\del P)
\subset
\del\usdoublestar(\del A_s\sm\del P)
$$
(cf. \cite{HN}, Theorem I).
Together with
Theorem \ref{RTThm} it follows that
over $A_s\sm\del P$ the function
$\usdoublestar$ is the solution of the HRMA (in dimension $n$)
$\h{\rm MA}\,\usdoublestar=0$. The difficulty in
replacing the regularity results of
\S\ref{ConeRegularitySection}--\ref{EssentialSmoothnessSection}
by the general $C^{1,1}$ regularity results for the Dirichlet
problem for the HRMA is that, aside from the fact that $A_s$ may be disconnected and
$P\sm A_s$ might not be convex,
it is not clear that $\del A_s$ will be
regular enough to apply the results of \cite{TU,CNS}.
Further, its boundary may intersect $\del P$, in which case
we would not be able to prescribe the Dirichlet data.

\bigs

\section{Partial $C^1$ regularity of the Legendre transform potential}
\label{PartialRegularityLegendrePotentialSection}

\bigs

Recall the definition of the Legendre transform potential,
$$
\psi(s,x)=\psi_s(x):= u_s^\star(x), \quad s\ge0,\quad x\in\RR^n.
$$
It is a one-parameter family of convex functions (in $x$),
and for each $s$ the function $\psi_s$ is defined and finite on all of
$\RR^n$ (see \S\ref{GradientSubDiffUsSection}).
Moreover, it is a convex function of $(s,x)$: by definition
\begin{equation}
\label{PsiSupDefEq}
\psi(s,x)
=
\sup_{y\in P}\big[\langle x,y\rangle-u_0(y)-s\dot u_0(y)\big]
=
\sup_{y\in P}\big[\big\langle\, (s,x),\,(-\dot u_0(y),y)\,\big\rangle-u_0(y)\big],
\end{equation}
i.e., $\psi$ is the supremum of linear functions in $(s,x)$,
hence convex.
Observe also that if we would have taken the Legendre transform of $u_0+s\dot u_0$
with respect to all $n+1$ variables the resulting function would be equal a.e. to
$+\infty$.

Set 
$$
\Sigma_\reg(T):=
\bigcup_{s\in[0,T]}\{s\}\times\del u_s\Big(\h{\rm int}\,\big(P\sm(\del P\cup A_s)\big)\Big)
\subseteq[0,T]\times\RR^n.
$$
Note that by Lemma \ref{LegendreLemma} we could have replaced the subdifferential
with the gradient in the definition of $\Sigma_\reg(T)$.

Given a convex function $f$, recall that
\begin{equation}
\label{DeltafDefEq}
\Delta(f)
\end{equation}
denotes the set on which $f$ is finite and differentiable.

The following result states that $\Sigma_\reg(T)$ 
coincides with the regular locus of $\psi$. Moreover,
it shows the following partial $C^1$ regularity for $\psi$---the 
regular set of $\psi$ is simply the (indexed) union 
of the regular sets of $\psi_s$.

\begin{prop}
We have
\label{SigmaRegCharProp}
$$
\Sigma_\reg(T)=\Delta_T(\psi):=\Delta(\psi)\;\cap\; [0,T]\times\RR^n.
$$
Further,
$$
\Delta(\psi)=\bigcup_{s>T_\span^\cvx}\{s\}\times\Delta(\psi_s).
$$

\end{prop}

Recall the definition of the singular locus of $\psi$ (\S\ref{IntroductionSection})
$$
\Sigma_\sing(T):=[0,T]\times\RR^n\;\sm (\Delta(\psi)\cap [0,T]\times\RR^n).
$$
The proposition gives the following explicit description
$$
\Sigma_\sing(T)=
\bigcup_{s\in[0,T]}\{s\}\times\nabla \usd(A_s\sm\del P)
\subseteq[0,T]\times\RR^n.
$$ 
Note that $\Delta(\psi)$ is not in general open (for instance,
consider the situation when a sequence of singular points 
$\{(s_k,x_k)\}_{k\ge1}\subset\Sigma_\sing(T)$ with 
$\lim_{k\ra\infty}s_k=T_\span^\cvx$ converges to
a point in $\Delta(\psi)$). 
However, it is dense, and its complement has Lebesgue
measure zero (\cite{Ro}, Theorem 25.5; \cite{RT}, Proposition 2.4).

Recall the following duality between differentiability 
and strict convexity.

\begin{lem}
\label{SmoothnessConvexityLemma}
{\rm (See \cite{Ro}, Theorem 26.3.)} A closed proper convex function is essentially
strictly convex if and only if its Legendre dual is essentially smooth.
\end{lem}

The proof of the Proposition will be a consequence of two lemmas proved below.

First, we show that $\psi$ is differentiable in $x$ on $\Sigma_\reg(T)$.
\begin{lem}
\label
{FirstConeRegularityPsiLemma}
 $(s,x)\in\Sigma_\reg(T)$ if and only if $x\in\Delta(\psi_s)$.
\end{lem}

\begin{proof}
By definition
$\psi_s(x)=\sup_{y'\in P}(\langle x,y'\rangle-\usdoublestar(y'))$
and the supremum is achieved when
$x\in\del \psi_s^\star(y')=\del\usdoublestar(y')=\nabla \usdoublestar(y')$.
Assume $(s,x)\in\Sigma_\reg(T)$. We claim that then
$y'\in\h{\rm int}\,\big(P\sm(\del P\cup A_s)\big)$.
Indeed, by our assumption we know that there exists
some $y\in\h{\rm int}\,\big(P\sm(\del P\cup A_s)\big)$
such that $x=\nabla u_s(y)$. Equation
(\ref{StrictConvexityUsdoublestarSupportSecondEq})
and Lemma \ref{LegendreLemma} (ii)
then imply that $\nabla \usdoublestar(y')\ne\nabla
\usdoublestar(y)$, unless $y'=y$, i.e.,
$y'=(\nabla\usd)^{-1}(x)$, as desired.
It follows that
\begin{equation}
\label
{LegendreDualFormulaRegPartPsiEq}
\psi_s(x)=\langle x,(\nabla\usd)^{-1}(x)\rangle-\usd\circ (\nabla\usd)^{-1}(x),
\quad x\in \Sigma_\reg(T).
\end{equation}
Since
$\del u_s\big(\h{\rm int}\,\big(P\sm(\del P\cup A_s)\big)\big)
=
\nabla u_s\big(\h{\rm int}\,\big(P\sm(\del P\cup A_s)\big)\big)$
is an open set, equation (\ref{LegendreDualFormulaRegPartPsiEq})
holds in a neighborhood of $x$ (for $s$ fixed).
Since $\usd$ equals $u_s$
on the open set $\h{\rm int}\,\big(P\sm(\del P\cup A_s)\big)$,
it is smooth there. Hence, we may differentiate
(\ref{LegendreDualFormulaRegPartPsiEq}) in $x$ to
obtain  (see (\ref{LegendreDualityGradientEq}))

\begin{equation}
\label{GradientInvertibleSigmaRegEq}
\nabla\psi_s(x)=(\nabla \usd)^{-1}(x),
\end{equation}
and this is a singleton.

Conversely, assume $x\in\Delta(\psi_s)$. Then
by Lemma \ref{SmoothnessConvexityLemma} it follows that
$\psi_s^\star$ is strictly convex in a neightborhood of $\nabla\psi_s(x)$.
We have
$\nabla\psi_s(x)\in \h{\rm int}\,\big( P\sm(\del P \cup A_s)\big)$
by Lemma \ref{WhereUsdStrictlyCvxLemma},
and 
$x\in
\del \usd\big(\h{\rm int}\,\big(P\sm(\del P\cup A_s)\big)\big)=
\del u_{s}\big(\h{\rm int}\,\big(P\sm(\del P\cup A_s)\big)\big)$
by duality (and Lemma \ref{LegendreLemma}),
i.e., $(s,x)\in\Sigma_\reg(T)$.
\end{proof}

Next, we show that wherever $\psi$ is differentiable in $x$ it is also
differentiable in $s$.

\begin{lem}
\label{PartialConeRegularityLemma}
Assume that $x\in\Delta(\psi_s)$. Then $(s,x)\in\Delta(\psi)$.
\end{lem}

\begin{proof}
It suffices to show that (\ref{LegendreDualFormulaRegPartPsiEq})
holds in a neighborhood of $(s,x)$.
Then the usual derivation
of the first variation formula for the Legendre transform
(see, e.g., \cite{R}, p. 85) gives
\begin{equation}
\label{FirstVariationLegendreSecondEq}
\frac{\del\psi}{\del s}(s,x)=-\frac{\del u}{\del s}(s,\nabla_x\psi(s,x))
=-\dot u_0\circ\nabla_x\psi(s,x),
\end{equation}
implying $(s,x)\in\Delta(\psi)$.

By the proof of Lemma \ref{FirstConeRegularityPsiLemma} we
already know that (\ref{LegendreDualFormulaRegPartPsiEq})
holds in a neighborhood of $x$, with $s$ fixed, and that,
further, it suffices to show that for some $\epsilon>0$ with $s+\epsilon<T$,
\begin{equation}
\label{SOpennessDelUsEq}
x\in\nabla u_{s'}\big(\h{\rm int}\,\big(P\sm(\del P\cup A_{s'})\big)\big),
\quad \forall\, s'\in(s-\epsilon,s+\epsilon),
\end{equation}
i.e., $(s',x)\in\Sigma_\reg(T)$.

First, we claim that if $y\in \h{\rm int}\,\big(P\sm(\del P\cup A_{s})\big)$
and $\epsilon>0$ is sufficiently small
then $y\in \h{\rm int}\,\big(P\sm(\del P\cup A_{s'})\big)$
for all $s'\in(s-\epsilon,s+\epsilon)$.
For that, it suffices to show that for $\epsilon>0$ sufficiently small
and $s'\in(s-\epsilon,s+\epsilon)$ there exists a neighborhood $U_{s,y}$
of $y$ such that
$\del u_{s'}(y')\ne\emptyset$ for all $y'\in U_{s,y}$.
Indeed, by (\ref{StrictConvexityUsdoublestarSupportFirstEq}),
\begin{equation}
\label{UsAboveSuppHyperplaneEq}
u_s(\tilde y)>\langle \tilde y-y,\nabla u_s(y)\rangle+u_s(y), \quad\forall
\tilde y\in P\sm\{y\},
\end{equation}
and moreover the same holds for $y$ replaced by a sufficiently nearby
point.
By smoothness of $\dot u_0$ and compactness of $P$ we also have that
$$
\max_{\tilde y\in P}|\dot u_0(\tilde y)-\dot u_0(y)-\langle \tilde y-y,\nabla\dot u_0(y)\rangle|
<C,
$$
with $C>0$ a uniform constant independent of $y$.
Hence
(\ref{UsAboveSuppHyperplaneEq}) holds with $s$ replaced by
$s'\in(s-\epsilon,s+\epsilon)$ and with $y$ replaced
by $y'$ with $|y-y'|<\epsilon$ for $\epsilon$ sufficiently small.
Thus, $\del u_{s'}(y')\ne\emptyset$, proving our claim.

Now, (\ref{SOpennessDelUsEq}) follows since
$x\in\nabla u_{s}\big(\h{\rm int}\,\big(P\sm(\del P\cup A_{s})\big)\big)$
and the mapping
$s\mapsto \nabla u_0+s\nabla\dot u_0$ is continuous in $s$.
\end{proof}

Combining the results obtained so far we may give
a proof of Proposition \ref{FirstMainProp}.

\medskip

\noindent
{\it Proof of Proposition \ref{FirstMainProp}.}
According to Lemma \ref{LegendreLemma} for each $s>0$ 
the function $\usd$ is finite, and according to 
Lemma \ref{SurjectivityRnLemma} its gradient surjects
to $\RR^n$. It follows that its dual $\psi_s$ must be
finite, since for every $x\in\RR^n$ the supremum in
the definition of the Legendre transform is necessarily
achieved for $y\in\del\psi_s(x)$ and any $y\in P$
satisfying $\del\usd(y)=x$ will do.
Convexity then implies that $\psi_s$ will be Lipschitz
continuous (see \cite{RT}, Proposition 2.4).

Lemma \ref{SmoothnessConvexityLemma}, combined with Lemma \ref{EssentialSmoothnessLemma},
Lemma \ref{WhereUsdStrictlyCvxLemma}, and the fact that $A_s\ne\emptyset$
for $s>T_\span^\cvx$, implies simultaneously that $\psi_s$ is 
essentially strictly convex,
and that it is not differentiable. 
In fact, since $\psi_s$ is convex and continuous on $\RR^n$
it follows that $\dom(\del \psi_s)=\RR^n$, thus
$\psi_s$ is strictly convex.
In addition, the exact description of the regular locus
of $\psi_s$ in Lemma \ref{FirstConeRegularityPsiLemma}
and the differentiability of $\usd$ (Lemma \ref{LegendreLemma})
imply that the singular locus of $\psi_s$ is
$\nabla\usd(A_s\sm\del P)$.
This concludes the proof of 
Proposition~\ref{FirstMainProp}.


\bigs

\section{The Legendre transform potential on the regular locus}
\label{LegendrePotentialRegularLocusSection}

\bigs

The following Proposition shows that $\psi$ solves the HRMA
wherever it is differentiable.
This is part of the statement of Theorem \ref{FirstMainThm}.

\begin{prop}
\label
{MAPsiRegPartProp}
(i)
For $T< T_\span^\cvx$,
\begin{equation}
\label{GraphSubDiffUptoTCvxEq}
\del\psi \big([0,T]\times\RR^n\big)
\;=\;
\h{\rm graph of $\;-\dot u_0\;$ over $\;P\sm\del P$}.
\end{equation}
(ii) Let $T>0$.
One has
\begin{equation}
\label{MAPsiRegularPartEq}
\h{\rm MA}\,\psi|_{\Sigma_\reg(T)}=0.
\end{equation}
Moreover,
\begin{equation}
\label{GraphSubDiffAfterTCvxRegEq}
\del\psi \big(\Sigma_\reg(T)\;\cap\; \{s\}\times\RR^n\big)
\;=\;
\h{\rm graph of $\;-\dot u_0\;$ over $\;\h{\rm int}\,\big(P\sm(\del P\cup A_s)\big)$}.
\end{equation}
(iii) $\Sigma_\reg(T)$ is dense in $[0,T]\times\RR^n$ and its complement has zero
Lebesgue measure there.
\end{prop}

\begin{proof}
(i)
For $T<T_\span^\cvx$ one has using  (\ref{VariationPotentialEq}),
\begin{equation}
\label{GraphDotUzeroPointsEq}
\begin{array}{rlr}
\del\psi(s,x)
& =
(
\dot\psi_s(x),
\nabla_x\psi_s(x)
)

=
(
-\dot u_s\circ\nabla_x\psi_s(x),
\nabla_x\psi_s(x)
)
\cr\cr
& =
(
-\dot u_0\circ\nabla_x\psi_s(x),
\nabla_x\psi_s(x)).
\end{array}
\end{equation}
Equation (\ref{GraphSubDiffUptoTCvxEq})
now follows from
$\nabla_x\psi_s(\RR^n)=P\sm\del P$.


\medskip
\noindent
(ii)
According to Lemmas \ref{FirstConeRegularityPsiLemma} and
\ref{PartialConeRegularityLemma} we know that $\psi$ is
differentiable on $\Sigma_\reg(T)$.
The same arguments as in the proof of Lemma \ref{FirstConeRegularityPsiLemma}
actually demonstrate that it is smooth on $\Sigma_\reg(T)$.
Differentiating (\ref{LegendreDualFormulaRegPartPsiEq}) twice
then yields (see (\ref{LegendreDualityHessianEq}))
$$
\big(\nabla^2\psi_s(x)\big)^{-1}=\nabla^2 u_s(\nabla\psi_s(x)),
$$
showing that $\psi_s$ is strictly convex at $x$ whenever
$(s,x)\in\Sigma_\reg(T)$.

Moreover, (\ref{LegendreDualFormulaRegPartPsiEq}) and
the proof of Lemma \ref{FirstConeRegularityPsiLemma}
show  that $\psi$ is in fact smooth on $\Sigma_\reg(T)$.
In particular, the second variation formula for the Legendre transform
\begin{equation}
\label{SecondVariationGeodesicFormulaEq}
-\ddot u_s|_{(\nabla u_s)^{-1}(x)}
=
\ddot\psi_s|_{x}+\langle\nabla\dot\psi_s|_{x},\nabla\dot u_s|_{(\nabla u_s)^{-1}(x)}\rangle
=
\big(\ddot\psi_s-{\textstyle\frac12}|\nabla\dot\psi_s|^2_{g_{\psi_s}}\big)(x),
\end{equation}
holds pointwise for $x\in\Sigma_\reg(T)$, and there these equations
are pointwise equivalent to the HRMA (\ref{HRMARayEq})
(see, e.g., \cite{R}, p. 87; \cite{RZ1}, \S3). Since $u_s$ solves the equation
on the left hand side, $\psi$ solves the HRMA on $\Sigma_\reg(T)$.

To prove (\ref{GraphSubDiffAfterTCvxRegEq}) observe that
(\ref{FirstVariationLegendreSecondEq}) holds on $\Sigma_\reg(T)$
(recall the arguments in Lemmas \ref{FirstConeRegularityPsiLemma} and
\ref{PartialConeRegularityLemma}).
Since, by duality,
$$
\del_x\psi\big(\Sigma_\reg(T)\;\cap\;\, \{s\}\times\RR^n\big)
=
\;\h{\rm int}\,\big(P\sm(\del P\cup A_s)\big),
$$
equation (\ref{GraphSubDiffAfterTCvxRegEq}) follows.
Note that (\ref{GraphSubDiffAfterTCvxRegEq}), together with
Theorem \ref{RTThm}, gives another proof
of  $\h{\rm MA}\,\psi\,|_{\Sigma_\reg(T)}=0$, since
the graph of $-\dot u_0$ over $P$ is a set of Lebesgue measure
zero in $\RR^{n+1}$.

\medskip
\noindent
(iii)
Observe that $\Sigma_{\sing}(T)$ is a set of Lebesgue measure zero
in $[0,T]\times\RR^n$ since
by Proposition \ref{SigmaRegCharProp}
it is precisely the set on which the convex
function $\psi$ is not differentiable.
\end{proof}

Observe that the proof of Proposition \ref{MAPsiRegPartProp} (iii)
relies on Lemma \ref{SurjectivityRnLemma} (i) implicitly.
Alternatively, one may also use Lemma \ref{SurjectivityRnLemma} (ii)
directly to prove that $\Sigma_\reg(T)$ is dense
in $[0,T]\times\RR^n$ without using the characterization
of $\Sigma_\reg(T)$ in Proposition \ref{SigmaRegCharProp}.

\bigs

\section{The subdifferential of the Legendre transform potential}
\label{SubdiffLegendrePotentialSection}

\bigs

In this section we prove upper and lower bounds on the 
total subdifferential of $\psi$
(Lemma \ref{YyPartialGammaPsiLemma}) in terms of 
$\usd$ and a partition of $A_s$.
To that end, we first 
collect some elementary facts concerning subdifferentials
of convex functions of several variables. Then, we 
give a precise description of the partial $x$-subdifferential
of $\psi$ and of the set of reachable partial $x$ subgradients
(Lemma \ref{YiPartialGammaPsiLemma}). This is then 
combined with the partial $C^1$ regularity of $\psi$ to
prove Lemma \ref{YyPartialGammaPsiLemma}.

Given a convex function $f$, define the set of reachable subgradients by
\begin{equation}
\label{GammaDefEq}
\gamma f(x):=\{x^*\,:\, \h{\rm \ exists\ }\, \{x_k\}_{k\ge 1}\subset\Delta(f) \h{\ \rm with\ }
\lim_{k\ra\infty}(x_k,\nabla f({x_k}))=(x, x^*)
\}.
\end{equation}
(recall (\ref{DeltafDefEq})).

\begin{lem}
\label{SubdifferentialConvexSetLemma}
{\rm (See \cite{HL1}, Theorem 6.3.1, p. 285.)}
Let $f$ be a closed proper convex function. Then
$\partial f(x)=\h{\rm co}\;\gamma f(x)$ for $x\in\dom(\del f)$.
\end{lem}

When $f(x_1,x_2)$ is convex in both of its arguments $x_1\in\RR^{m_1}, x_2\in\RR^{m_2}$,
we will denote the reachable partial $x_1$ subgradients by
$$
\gamma_{x_1} f(x_1,x_2):=\gamma f_{x_2}(x_1),
$$
where $f_{x_2}(x_1):=f(x_1,x_2)$ is considered as a function of $x_1$.
Similarly we also define $\gamma_{x_2}f(x_1,x_2)$.
We will denote by $\del_{x_1}f$ the partial $x_1$ subdifferential of
$f$ considered as a function of $x_1$ alone,
$$
\del_{x_1}f(x_1,x_2):=\del f_{x_2}(x_1).
$$
Similarly we also define $\del_{x_2}f$.

We will need the following elementary consequence of Lemma \ref{SubdifferentialConvexSetLemma}.

\begin{lem}
\label{TotalSubdifferentialProp}
Let $f(x_1,x_2)$, $x_1\in\RR^{m_1},x_2\in\RR^{m_2}$, be a closed proper
convex
function on
$\RR^{m_1+m_2}$. Assume that $f$ is sub-differentiable, and differentiable in $x_1$. Then
wherever $f$ is finite
\begin{equation}
\label{SubdifferentialFEq}
\partial f = \nabla_{x_1} f\times\partial_{x_2} f.
\end{equation}
\end{lem}

\begin{proof}
Fix $x=(x_1,x_2)$.
By Lemma \ref{SubdifferentialConvexSetLemma}, we have
\begin{equation}
\label{SubdifferentialPsiContainmentEq}
\begin{array}{lll}
\partial f{(x_1,x_2)}
& \supset &
\{(x_1^*,x_2^*)\,:\, x_1^*\in\gamma_{x_1}f{(x_1,x_2)},
x_2^*\in \gamma_{x_2} f{(x_1,x_2)}\}
\cr\cr
& = &
\nabla_{x_1}f(x)\times\gamma_{x_2} f{(x_1,x_2)},
\end{array}
\end{equation}
Indeed one takes sequences of points where $f$ is
differentiable and for these points $\nabla f=(\nabla_{x_1}f,\nabla_{x_2}f)$.
This is possible since a closed proper function is differentiable on a dense
set in the interior of $\dom(f)$, the set where it is sub-differentiable
(\cite{Ro}, Theorem 25.5).
Lemma \ref{SubdifferentialConvexSetLemma} implies that $\partial f(x)$
is a convex set, and by applying this Lemma yet once more,
that is taking convex hulls in (\ref{SubdifferentialPsiContainmentEq}),
we obtain
$$
\partial f{(x_1,x_2)}
\supset
\nabla_{x_1}f(x)\times\partial_{x_2} f(x).
$$
From the definitions one may verify that
\begin{equation}
\label{SubDiffProductEq}
\partial f{(x_1,x_2)}
\subset
\del_{x_1}f(x)\times\partial_{x_2} f(x)
=
\nabla_{x_1}f(x)\times\partial_{x_2} f(x),
\end{equation}
and this completes the proof.
\end{proof}

The following lemma describes a general property of convex functions of several variables,
that we prove using the concept of reachable subgradients.
It also appears, with a different proof using the Hahn-Banach theorem, 
in \cite{A}, Proposition 2.4.

\begin{lem}
\label
{SurjectiveProjectionLemma}
The projection $\pi_x:\del\psi(s,x)\mapsto\del_x\psi(s,x)$
is surjective.
\end{lem}

\begin{proof}
First, as in (\ref{SubDiffProductEq}), we have

\begin{equation}
\label{SubDiffProductSecondEq}
\del\psi(s,x)\subset\del_x\psi(s,x)\times\del_s\psi(s,x),
\end{equation}
and so $\pi_x$ indeed maps $\del\psi(s,x)$ into $\del_x\psi(s,x)$.
The set $\del\psi(s,x)\subset\RR^{n+1}$ is convex.
Therefore, by Lemma \ref{SubdifferentialConvexSetLemma},
it suffices to verify that the projection $\pi_x$
surjects to $\gamma_x\psi(s,x)$.

Let $y\in\gamma_x\psi(s,x)$. Then there exists a sequence of points
$\{x_i\}_{i\ge1}$ converging to $x$ and such that $\nabla_x\psi(s,x_i)$
exists and $\lim_{i\ra\infty}\nabla_x\psi(s,x_i)=y$.
By Lemma \ref{TotalSubdifferentialProp}
$$
\lim_{i\ra\infty}\del\psi(s,x_i)
=
\lim_{i\ra\infty}\;\;\del_s\psi(s,x_i)\times \{\nabla_x\psi(s,x_i)\},
$$
from which $\pi_x\big(\lim_{i\ra\infty}\del\psi(s,x_i)\mskip1mu\big)=y$.
By the upper semi-continuity of the subdifferential mapping (\cite{Ro}, Corollary 24.5.1)
then $\lim_{i\ra\infty}\del\psi(s,x_i)\subset\del\psi(s,x)$.
We conclude that $y\in\pi_x\big( \del\psi(s,x) \big)$, as claimed.
\end{proof}

For each $(s,y)\in\RR_+\times A_s\sm\del P$, set
\begin{equation}
\label
{QMaxSYDefEq}
Q(s,y):=\{v\in P\,:\, \nabla\usdoublestar(v)=\nabla\usdoublestar(y)\}
\subset \overline{A_s}\sm\del P
\end{equation}
(the inclusion is implied by Lemma \ref{EssentialSmoothnessLemma} 
and the results of \S\ref{InvertibilitySection}),
and
\begin{equation}
\label{YyEq}
Y(s,y):=Q(s,y)\cap (\del A_s\sm\del P)=Q(s,y)\cap \del A_s.
\end{equation}
Note that $Q(s,y)$ is the projection of the intersection of
the tangent hyperplane to
the graph of $\usdoublestar$ at $y$ with the graph itself.
Thus, convexity of $\usdoublestar$ implies that $Q(s,y)$ is closed and convex, and
\begin{equation}
\label{QMaxCoYsyEq}
Q(s,y)=\h{\rm co}\, Y(s,y)\subset P
\end{equation}
(here we use that any point in $\h{\rm int}\,(A_s\sm\del P)$
is a convex combination of called upon points lying in $\del A_s$,
see the proof of Lemma \ref{LegendreLemma}).

Note also that $\{y_i\}_{i\in I}\subset Y(s,y)$.
Recall that the the set $\{y_i\}_{i\in I}$ of points
called upon by $y$ obtained in the proof of
Lemma \ref{LegendreLemma} was not necessarily the unique
collection of at most $n+1$ points satisfying
(\ref{CalledUponEq}) and (\ref{ISetEq}).
The set $Y(s,y)$ will serve as a more `canonical'
substitute for $\{y_i\}_{i\in I}$.

\begin{lem}
\label
{YiPartialGammaPsiLemma}
Let $y\in A_s\sm\del P$, let $Y(s,y)$ be given by (\ref{YyEq}),
and let $x=\nabla \usdoublestar(y)$. Then
\begin{equation}
\label{YiPartialGammaPsiEq}
Y(s,y)=\gamma_x\psi(s,x),
\end{equation}
and
\begin{equation}
\label{QMaxPartialDelPsiEq}
Q(s,y)=\del_x\psi(s,x).
\end{equation}

\end{lem}

\begin{proof}
Let $v\in Y(s,y)$.
Since $v\in\del A_s\sm\del P$ (Claim \ref{CalledUponInteriorClaim})
there exists
a sequence
$\{v_k\}_{k\ge1}\subset \h{\rm int}\,\big( P\sm(\del P\cup A_s)\big)$
with $\lim_{k\ra\infty}v_k=v$. By (\ref{GradientInvertibleSigmaRegEq})
$\nabla\psi_s$ exists and is invertible on
$\h{\rm int}\,\big( P\sm(\del P\cup A_s)\big)$.
So for each $k$ there exists a unique $x_k$
with
\begin{equation}
\label{VkXkEq}
\nabla u_s(v_k)=(\nabla\psi_s)^{-1}(v_k)=x_k.
\end{equation}
Therefore, $\lim_{k\ra\infty} x_k=\nabla u_s(v)=x$, and
$\lim_{k\ra\infty}\nabla\psi_s(x_k)=v$.
It follows that
$Y(s,y)\subset\gamma_x\psi(s,x)$.

Conversely, if $v\in\gamma_x\psi(s,x)$ then
by definition there exists a sequence
$\{x_k\}_{k\ge1}$ in $\Delta(\psi_s)$ with
$\lim_{k\ra\infty}\nabla\psi_s(x_k)=v$. By Lemma
\ref{FirstConeRegularityPsiLemma} $(s,x_k)\in\Sigma_\reg(T)$
for all $k$. By duality then
$\nabla\psi_s(x_k)\in \h{\rm int}\,\big( P\sm (\del P\cup A_s)\big)$.
Since $(s,x)\in\Sigma_\sing(T)$ we must have
$\del_x\psi(s,x)\ni v\not\in \h{\rm int}\,\big( P\sm (\del P\cup A_s)\big)$.
We conclude that $v\in\del A_s$.
Finally, $\del_x\psi(s,x)\subset Q(s,y)$:
if $v\in\del_x\psi(s,x)$ then by duality, $x\in\del\usdoublestar(v)$ and by
Lemma \ref{LegendreLemma} $x=\nabla\usdoublestar(v)$.
Therefore, $v\in Y(s,y)$, implying (\ref{YiPartialGammaPsiEq}).

Taking convex hulls in (\ref{YiPartialGammaPsiEq}) and invoking
(\ref{QMaxCoYsyEq}) and Lemma \ref{SubdifferentialConvexSetLemma}
implies (\ref{QMaxPartialDelPsiEq}).
\end{proof}

\begin{lem}
\label
{YyPartialGammaPsiLemma}

Let $y\in A_s\sm\del P$, 
and set $x:=\nabla\usdoublestar(y)$.
Then
\begin{equation}
\label{TildeQMaxContainedDelPsiEq}
\h{\rm co}\,\big\{(-\dot u_0(v),v)\,:\, v\in \gamma_x\psi(s,x)\big\}
\subset
\del\psi(s,x)
\subset
\h{\rm co}\,\big\{(-\dot u_0(v),v)\,:\, v\in\del_x\psi(s,x)\big\}.
\end{equation}

\end{lem}

\begin{proof}
First, we will show that
\begin{equation}
\label{GammaPsiEqualYsyEq}
\big\{(-\dot u_0(v),v)\,:\, v\in \gamma_x\psi(s,x)\big\}
\subset
\gamma\psi(s,x).
\end{equation}
Together with Lemma \ref{SubdifferentialConvexSetLemma}
this implies the first inclusion in (\ref{TildeQMaxContainedDelPsiEq}).
For that, let $v\in Y(s,y)$.
Let $\{v_k\}_{k\ge1}$ and $\{x_k\}_{k\ge1}$  be as
in the proof of Lemma \ref{YiPartialGammaPsiLemma}.
By Lemma \ref{PartialConeRegularityLemma} we have
$(s,x_k)\in\Delta(\psi)$
and $$
\nabla\psi(s,x_k)=(-\dot u_0\circ\nabla\psi_s(x_k),\nabla\psi_s(x_k)).
$$
Hence, by (\ref{VkXkEq}) 
$\lim_{k\ra\infty}\nabla\psi(s,x_k)=(-\dot u_0(v),v)\in\gamma\psi(s,x)$,
and (\ref{GammaPsiEqualYsyEq}) follows.

Next, let $\tilde v\in\gamma\psi(s,x)\subset\del\psi(s,x)$, and
let $\{(s_k,x_k)\}_{k\ge1}\subset \Delta(\psi)$
satisfy
$$
\lim_{k\ra\infty}(s_k,x_k)=(s,x),\quad
\lim_{k\ra\infty}\nabla\psi(s_k,x_k)=\tilde v.
$$
Then, as in (\ref{LegendreDualityEq})-(\ref{LegendreDualityGradientEq}) (since
we are in $\Delta(\psi)$),
\begin{equation}
\label{NablaPsiSeqKEq}
\nabla\psi(s_k,x_k)
=
\big(\dot \psi(s_k,x_k),\nabla_x\psi({s_k},x_k)\big)
=
\big(-\dot u_0\circ\nabla_x\psi_{s_k}(x_k),\nabla_x\psi_{s_k}(x_k)\big).
\end{equation}
Let $v:=\lim_{k\ra\infty}\nabla_x\psi(s_k,x_k)\in\del_x\psi(s,x)$
(making use of (\ref{SubDiffProductSecondEq})).
Then (\ref{NablaPsiSeqKEq}) implies that
$$
\tilde v=(-\dot u_0(v),v).
$$
The second inclusion in (\ref{TildeQMaxContainedDelPsiEq})
now follows from Lemma \ref{SubdifferentialConvexSetLemma}.
\end{proof}

We note that an alternative proof, starting from
the formula (\ref{PsiSupDefEq}), may be given.
Also, the lower bound can in fact be shown to be sharp. We do not
go into the details here since we will not need these facts
in the present article.

\bigs
\section{Monge-Amp\`ere measure of the Legendre transform potential}
\label{MAMeasureLegendrePotentialSection}

\bigs
We are now in a position to complete the proof of Theorem \ref{FirstMainThm},
and show that the Legendre transform potential $\psi$
is not a weak solution of the HRMA (\ref{HRMARayEq}) for $T>T_\span^\cvx$.

Theorem \ref{FirstMainThm} is a direct corollary of the following result,
stating that the \MA measure of $\psi$ charges the singular locus of $\psi$
with positive mass.

\begin{prop}
\label{MAMassPsiProp}
Let $T>T_\span^\cvx$. Then
\begin{equation}
\label{MAMassEq}
\int_{[0,T]\times\RR^n}\h{\rm MA}\,\psi
=
\int_{\Sigma_\sing(T)}\h{\rm MA}\,\psi
>0.
\end{equation}
\end{prop}

The proof of Proposition \ref{MAMassPsiProp} will occupy most of this
section.

Set
\begin{equation}
\label{TildeQMaxSYDefEq}
\tilde Q(s,y):=
\h{\rm co}\,\big\{(-\dot u_0(v),v)\,:\, v\in Y(s,y)\big\}
\subset \RR\times P.
\end{equation}
Note that by Lemmas \ref{YiPartialGammaPsiLemma} and
\ref{YyPartialGammaPsiLemma} this convex set is contained
in $\del\psi(s,\nabla\usdoublestar(y))$.

Let
$$
q(s,y):=\dim Q(s,y)\le n.
$$
Recall that
$$
\{y_i\}_{i\in I}\subset Y(s,y)\subset Q(s,y),
$$
and that $|I|>1$, whenever $y\in A_s$. Hence, when $y\in A_s$,
the set $Q(s,y)$ contains at least a line,
and $q(s,y)\ge 1$.
By (\ref{QMaxCoYsyEq}) we have $\pi_P(\tilde Q(s,y))=Q(s,y)$.
Hence, $q(s,y)\le\dim\tilde Q(s,y)\le q(s,y)+1$.

\begin{lem}
\label{EpigraphdotuzeroLemma}
Let $y\in A_s\sm\del P$.
Then
\begin{equation}
\label
{PiPProjectedSetEq}
U(s,y):=\pi_P
\Big(
\tilde Q(s,y)
\;\sm\;
\big\{
(-\dot u_0(v),v)\,:\, v\in Q(s,y)
\big\}
\Big)\subset Q(s,y)
\end{equation}
is a set with a non-empty interior relative to $Q(s,y)$.
If
 $\dim\tilde Q(s,y)=q(s,y)$ then
\begin{equation}
\label
{UsyIntersectionEq}
\tilde Q(s,y)
\cap
\{(-\dot u_0(v),v)\,:\, v\in U(s,y)\}
=
\emptyset.
\end{equation}

\end{lem}

\begin{proof}
We may assume that $\dim\tilde Q(s,y)=q(s,y)$ since
if $\dim\tilde Q(s,y)=q(s,y)+1$ then in fact the set in
(\ref{PiPProjectedSetEq}) contains $\h{int}\,Q(s,y)$.

By our assumption $\tilde Q(s,y)$ contains at most $q(s,y)+1$
linearly independent points in $\RR^{n+1}$.
By convexity, $\tilde Q(s,y)$ must be
an affine graph over $Q(s,y)$.

Let $v_0,v_1\in Q(s,y)\cap \del A_s$ (necessary in $\del A_s\sm\del P$ by
Lemma \ref{EssentialSmoothnessLemma}) and assume $v_0\ne v_1$.
This is possible since $\{y_i\}_{i\in
I}\subset Q(s,y)\cap \del A_s$
and $|I|>1$.
Consider the function
$$
F(a):=u_s\big((1-a)v_0+av_1\big) \quad a\in[0,1].
$$
Denote $v(a):=(1-a)v_0+av_1\in Q(s,y)$ (the inclusion follows
from convexity of $Q(s,y)$). Then
\begin{equation}
\label
{FijSecondDerivativeEq}
F''(a)
=
\langle \nabla^2 u_0|_{v(a)}.(v_1-v_0),v_1-v_0\rangle
+
s\langle \nabla^2 \dot u_0|_{v(a)}.(v_1-v_0),v_1-v_0\rangle.
\end{equation}
If the statement of the Lemma were false then in fact by
continuity it would necessarily follow that
$$
\tilde Q(s,y)=\h{\rm graph of \ $-\dot u_0$\ \ over\ \ $Q(s,y)$}.
$$
However, since $\tilde Q(s,y)$ is affine over $Q(s,y)$
then the second term in (\ref{FijSecondDerivativeEq}) would
vanish, implying strict convexity of $F$ on $[0,1]$. However, that
contradicts the fact
that $\nabla u_s(v_0)=\nabla u_s(v_1)$
(see (\ref{GradientSubDiffBasicEq})
and (\ref{QMaxSYDefEq})), hence $F'(0)=F'(1)$.

Finally, (\ref{UsyIntersectionEq}) follows from (\ref{PiPProjectedSetEq})
since, if $\dim\tilde Q(s,y)=q(s,y)$,
the fibers of the projection $\pi_P:\tilde Q(s,y)\ra Q(s,y)$
are singletons.
\end{proof}


\vglue1cm

\hglue2.26cm\includegraphics[scale=0.5]{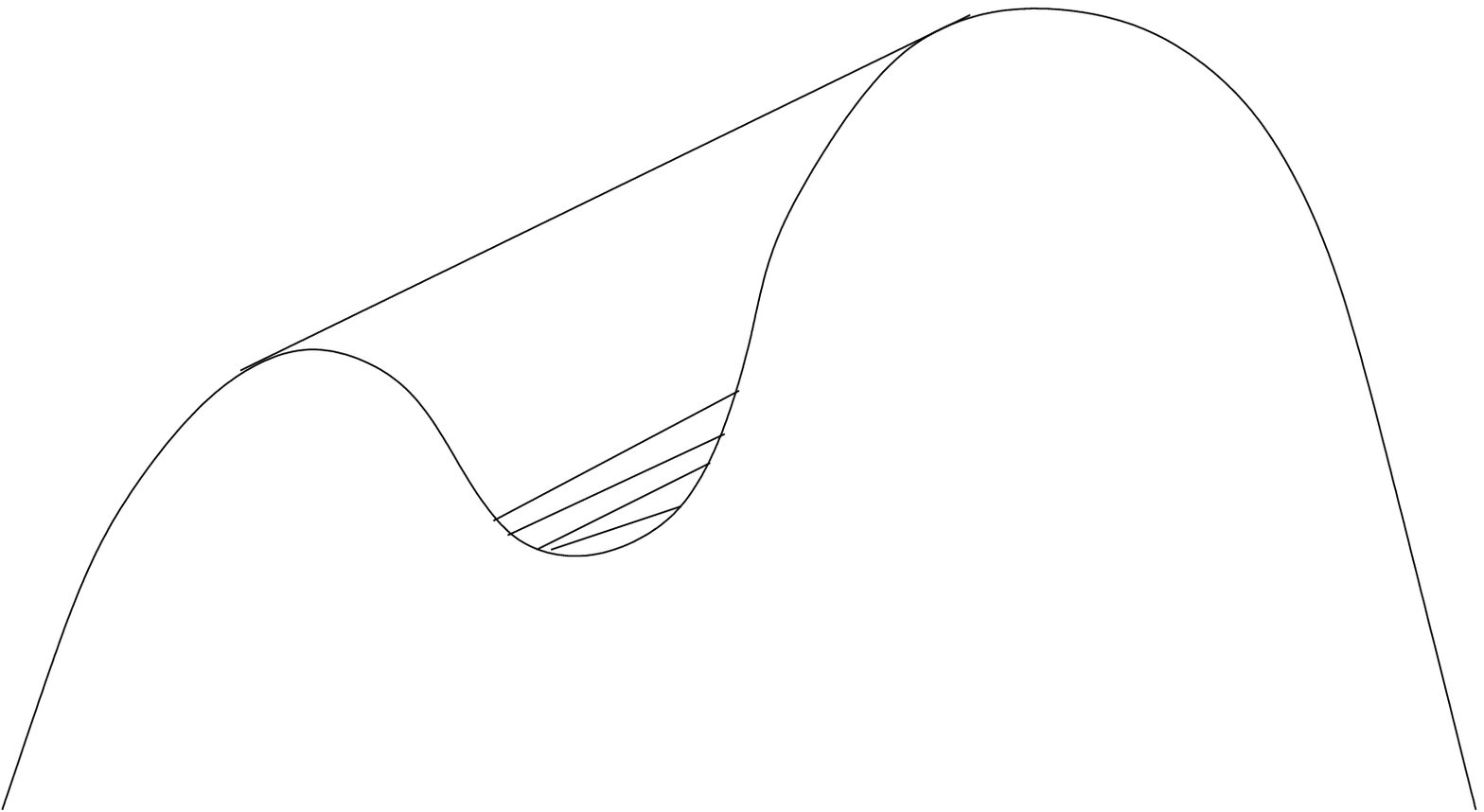}


\bigskip\bigskip
\noindent
{\sc Figure 4.} The graphs of $-\dot u_0$ and $-(\dot u_0^{\star\star})$
over $P$, and line segments in the image
\hfill\break
 of $\del\psi$ corresponding to $\tilde Q(s,y)$ for different values
of $s$.

\vglue1.5cm

\begin{lem}
\label
{MASingMassAinfinityLemma}
Let $T>T_\span^\cvx$. Then

\bigskip
\noindent
\begin{equation}
\label{MAMassSecondEq}
0<\int_{\Sigma_\sing(T)}\h{\rm MA}\,\psi
<\h{\rm Vol}\, 
\big(\,
 \epi\,(-\dot u_0)\, \sm \epi\,(-(\dot u_0)^{\star\star})
\,\big).
\end{equation}
\end{lem}

\begin{proof}
Let $s>T_\span^\cvx$ and let $y\in A_s\sm\del P$.
Note that we have a partition
\begin{equation}
\label{AsPartitionEq}
\overline{A_s}\sm\del P=\bigcup_{v\in A_s\sm\del P}Q(s,v),
\end{equation}
that is any two sets appearing in the union are
either disjoint or else coincide.
Let $x:=\nabla\usdoublestar(Q(s,y))$.
By Lemma \ref{YiPartialGammaPsiLemma}
$Q(s,y)=\del_x\psi(s,x)$, and further by
Lemma \ref{YyPartialGammaPsiLemma},
\begin{equation}
\label{TildeQMaxContainedDelPsiSecEq}
\tilde Q(s,y)=
\h{\rm co}\,\big\{(-\dot u_0(v),v)\,:\, v\in Y(s,y)\big\}
\subset
\del\psi(s,x).
\end{equation}
Therefore, from (\ref{AsPartitionEq}) we conclude that
\begin{equation}
\label{SetContainedInDelPsiEq}
\del\psi\big(\{s\}\times\RR^n\big)
\supset
\bigcup_{v\in A_s\sm\del P}\tilde Q(s,v).
\end{equation}

According to Proposition \ref{MAPsiRegPartProp} (i) the set
$\del\psi([0,T_\span^\cvx)\times\RR^n)$ is equal to the graph
of $-\dot u_0$ over $P\sm\del P$. On the other hand,
whenever $T>T_\span^\cvx$, the set on the right hand side
of (\ref{SetContainedInDelPsiEq})
projects to the non-empty set $\overline{A_s}\sm\del P$.
Above each piece $Q(s,v)$ of the partition of $A_s\sm\del P$,
the set $\tilde Q(s,v)$ satisfies
$\pi_P(\tilde Q(s,v))=Q(s,v)$
and is either a convex body of one dimension higher than $Q(s,v)$
or else is an affine graph over $Q(s,v)$.
We now describe these two possibilities in more detail.

Let $s\in(T_\span^\cvx,T)$ and let $y\in A_s\sm\del P$.
Either,

\smallskip\noindent
(a) $\dim\tilde Q(s,y)=q(s,y)$, in which case
by Lemma \ref{EpigraphdotuzeroLemma} there
exists a subset $U(s,y)$ of $Q(s,y)$, with non-empty
interior relative to $Q(s,y)$, such that
$$
\tilde U(s,y):=
\tilde Q(s,y)\cap \pi_P^{-1}(U(s,y))
$$
is affine and does not intersect the graph of $-\dot u_0$
over $U(s,y)$.
Or,

\smallskip\noindent
(b) $\dim\tilde Q(s,y)=q(s,y)+1$ (while the graph of $-\dot u_0$ over
$Q(s,y)$ has dimension $q(s,y)$). In this 
case there exist $q(s,y)+1$ distinct points 
$\big\{v_i\big\}_{i=1}^{q(s,y)+1}\subset Y(s,y)$ such that
$Q'(s,y):=\h{co}\,\big\{v_i\big\}_{i=1}^{q(s,y)+1}\subset Q(s,y)$ is a polyhedron
of dimension $q(s,y)$, and such that
$\big\{(-\dot u_0(v_i),v_i)\big\}_{i=1}^{q(s,y)+1}\subset\tilde Q(s,y)$,
are linearly independent in $\RR^{n+1}$.
Then $\tilde Q'(s,y):=\h{co}\,\big\{(-\dot u_0(v_i),v_i)\big\}_{i=1}^{q(s,y)+1}\subset\tilde Q(s,y)$
is a polyhedron of dimension $q(s,y)$ that can be considered
as an affine graph over $Q'(s,y)$. Imitating the proof of 
Lemma \ref{EpigraphdotuzeroLemma} for $Q'(s,y)$ instead of 
$Q(s,y)$ then implies that there exists 
a set $U(s,y)\subset Q'(s,y)$ with non-empty interior relative
to $Q'(s,y)$, and hence also relative to $Q(s,y)$, such that
if we set
$$
\tilde U(s,y):=
\tilde Q'(s,y)\cap \pi_P^{-1}(U(s,y))
$$
then
$$
\tilde U(s,y)\cap \big\{(-\dot u_0(v),v)\,:\, v\in U(s,y)\big\}=\emptyset.
$$

Note that the construction in (b) in effect 
reduces case (b) to case (a).

\medskip\noindent
Let
$$
V(s):=\bigcup_{v\in A_s}U(s,v)\subset A_s\sm\del P.
$$
It follows from the above 
and a Fubini type theorem
that for each $s>T_\span^\cvx$ the non-empty set $V(s)$
has positive Lebesgue measure in $\RR^n$ (note that,
in principle, as defined, $V(s)$ might not be open): indeed, each
$U(s,y)$ is open relative to $Q(s,y)$ and by (\ref{AsPartitionEq})
the union of the sets $Q(s,v)$ over $v\in A_s\sm\del P$ 
equals $\overline{A_s}\sm\del P$, that 
has non-empty interior in $\RR^n$. 
We set
$$
\tilde V(s):=
\bigcup_{v\in A_s}\tilde U(s,y)\subset \RR\times A_s\sm\del P.
$$
By construction $\tilde V(s)$ is a locally affine graph
over $V(s)$.

\begin{definition}
For each $s>T_\span^\cvx$, let $f(s,\,\cdot\,):A_s\ra\RR$ 
denote the unique locally affine function whose graph
over $V(s)$ equals $\tilde V(s)$ and that is affine
over each of the sets $Q(s,y)$.

\end{definition}

The graph of $f(s,\,\cdot\,)$ is obtained by 
extending affinely each affine piece $\tilde U(s,y)$ originally
defined over $U(s,y)$ to all of $Q(s,y)$.
Note that the graph of $f(s,\,\cdot\,)$ restricted to $V(s)$ satisfies 
$$
\emptyset
=
\big\{(-\dot u_0(v),v)\,:\, v\in V(s)\big\}
\;\cap\;
\big\{(f(s,v),v)\,:\, v\in V(s)\big\}.
$$
Further, $\del\psi(\{s\}\times\RR^n)$
contains $\tilde V(s)$, i.e., the graph of $f(s,\,\cdot\,)$ restricted to $V(s)$.

Now, by Lemmas \ref{AinfinityLemma} and \ref{AsContinuousMappingLemma}
the sets $A_s$ vary continuously and are strictly monotonically increasing
in $s$. Also, $A_{s'}=\emptyset$ for $s'<T_\span^\cvx$.
In addition, each piece $Q(s,v)$ intersects $\del A_s$. 
It follows that we may assume that
$f(s,y)$ defined above is continuous in $s$ for a.e. $y\in A_s$.
Furthermore, it also follows that there exists a set
$I\subset (T_\span^\cvx,T]$ of positive measure (in $\RR$)
and subsets $V'(s)\subset V(s)$ of positive
Lebesgue measure (in $\RR^n$), for each $s\in I$,
such that
$$
\big\{(f(s,v),v)\,:\, v\in V'(s)\big\}
\;\cap\;
\big\{(f(s',v),v)\,:\, v\in V'(s')\big\}
=\emptyset,
\quad s\not=s',\quad s,s'\in I.
$$
Let $F(s,y)dy$ be the volume measure induced on the graph of
$f(s,y)$ over $A_s$, regarded as a hypersurface in $\RR^{n+1}$,
from the Lebesgue measure on $\RR^{n+1}$. From the above
we have $\int_{V'(s)}F(s,y)dy>0$ for each $s>T_\span^\cvx$.
Fubini's theorem now gives that
$$
\int_{T_\span^\cvx}^Tds\int_{V(s)}F(s,y)dy
>
\int_{I}ds\int_{V'(s)}F(s,y)dy
>0.
$$
Hence,
$$
\bigcup_{{s\in(T_\span^\cvx,T]}\atop{v\in V(s)}}
\big\{(f(s,v),v)\,:\, v\in V(s)\big\}
$$
and therefore also
$$
\bigcup_{{s\in(T_\span^\cvx,T]}\atop{v\in A_s\sm\del P}}
\tilde Q(s,v)
$$
contains a set of positive Lebesgue measure in $\RR^{n+1}$.
By (\ref{SetContainedInDelPsiEq}) this set is
contained in $\del\psi\big((T_\span^\cvx,T]\times\RR^n\big)$.
According to Theorem \ref{RTThm} this implies the lower bound in
(\ref{MAMassSecondEq}).

The upper bound in (\ref{MAMassSecondEq})
follows from Lemma \ref{AinfinityLemma} and the upper bound
of Lemma \ref{YyPartialGammaPsiLemma}.
\end{proof}

Lemma \ref{MASingMassAinfinityLemma}, together with
Proposition \ref{MAPsiRegPartProp}  conclude the proof
of Proposition \ref{MAMassPsiProp},
from which Theorem \ref{FirstMainThm} follows.

\bigskip

\noindent {\bf Acknowledgments.}
This material is based upon work supported in part under National Science Foundation
grants DMS-0603850, 0904252.
Y.A.R. was also supported by a National Science Foundation
Postdoctoral Research Fellowship at
Johns Hopkins University during the academic year 2008--2009.



\end{document}